\documentclass[11pt]{article}
\usepackage[latin1]{inputenc}
\usepackage[dvips,final]{graphics}
\usepackage[french]{babel}
\usepackage{amssymb}
\usepackage{amsfonts}
\usepackage{amstext}
\usepackage{amsopn}
\usepackage{fancybox,color,epsf}
\usepackage{rotating}

\catcode`\ =\active \def {\`e} \catcode`\ =\active \def {\`a} \catcode`\ =\active \def {\`u} \catcode`\
=\active \def {\'e} \catcode`\ =\active \def {\c} \def \parn{\par\noindent} \catcode`\ =\active \def {\^e}
\catcode`\ =\active \def {\^o} \catcode`\ =\active \def {\^\i } \catcode`\ =\active \def {\^u}

\catcode`\: =\active \def :{\ifmmode\string:\else\unskip\kern 2pt\string: \ignorespaces\fi } \catcode`\; =\active
\def ;{\ifmmode\string;\else\unskip\kern 2pt\string; \ignorespaces\fi }

\def \m{\mu}  \def \n{\nu}
  \def \b{\beta}     \def \d{\delta}     
\def \la{\lambda} \def \La{\Lambda}  \def \s{\sigma}   \def \t{\theta}  
\def \e{\varepsilon}       

\def \x{\mathbf{x}} \def \p{\mathbf{p}} \def \v{\mathbf{v}}
 \def \W{\mathbf{W}} 
\def \M{\mathbf{M}}

\def \E{\mathbb{E}}
\def \P{\mathbb{P}}

\def \pa{\partial}        \def \ii{\infty}    
          \def \Ra{\Rightarrow}
       
              \def \nea{\nearrow}        \def \sea{\searrow}

\def \ub{\underbar}      

   \def \ds{\displaystyle}   \def \ts{\textstyle}       
\def \rt1{\sqrt{-1}\,\,}  \def \1{^{-1}}            \def \2{^{-2}}             \def \5{{\ts {1\over 2}}}
\def \parn{\par\noindent}

\begin{document}

\newtheorem{defi}{Definition}
\newtheorem{theo}{Theorem}
\newtheorem{prop}{Proposition}
\newtheorem{propr}{Property}
\newtheorem{cor}{Corollary}
\newtheorem{lem}{Lemma}
\newtheorem{rem}{Remark}

\hspace{1cm} \par
\if{
\title{\bf Central Limit Theorem for a class of Relativistic Diffusions}
\author{J\"urgen ANGST \  and \  Jacques FRANCHI}
\maketitle 
}\fi 

\centerline{{\bf{\LARGE{Central Limit Theorem for a Class of Relativistic Diffusions} }}} \par 
\vspace{0.6cm} \centerline{\textsc{ {J\"urgen ANGST \  and \  Jacques FRANCHI}}} \par

\if{
\centerline{{\textsc{\LARGE{central limit theorem for a class} }}} \par
\vspace{0.1cm} \centerline{\textsc{\LARGE{of relativistic diffusions}}} \par
\vspace{0.5cm} \centerline{\textsc{ {J\"urgen ANGST \  and \  Jacques FRANCHI}}} \par
}\fi 
\vspace{0.8cm}
\centerline{
\begin{minipage}{13cm}
\begin{sl}\small{\bf Abstract }:
Two similar Minkowskian diffusions have been considered, on one hand by Barbachoux, Debbasch, Malik and Rivet ([BDR1], [BDR2], [BDR3], [DMR], [DR]), and on the other hand by Dunkel and Hänggi ([DH1], [DH2]). We address here two questions, asked in [DR] and in ([DH1], [DH2]) respectively, about  the asymptotic behaviour of such diffusions. More generally, we establish a central limit theorem for a class of Minkowskian diffusions, to which  the two above ones belong. As a consequence, we correct a partially wrong guess in [DH1]. 
\end{sl}
\end{minipage}}


\section{Introduction} 

    Debbasch, Malik and Rivet introduced in [DMR] a relativistic diffusion in Minkowski space, they called Relativistic Ornstein-Uhlenbeck Process (ROUP), to describe the motion of a point particle surrounded by a heat bath, or relativistic fluid, with respect to the rest-frame of the fluid, in which the particle diffuses.  This ROUP was then studied in [BDR1], [BDR2], [BDR3], [DR], and extended to the curved case in [D]. Then Dunkel and Hänggi introduced and discussed in [DH1], [DH2] a similar process, in Minkowski space, they called relativistic Brownian Motion.  \par 
    Note that independently, a Relativistic Diffusion on any Lorentz manifold was defined in [FLJ], as the only diffusion whose law possesses the relativistic invariance under the whole isometry group of the manifold.  Accordingly, the particular case of the Schwarzschild-Kruskal-Szekeres manifold was studied. The case of G\" odel's universe was recently studied in [F]. \par 
   In [DR], Debbasch and Rivet argue qualitatively that the so-called ``hydrodynamical limit'' of their ROUP should behave in a Brownian way. They stress that a mathematical rigourous proof remains needed, to confirm such not much intuitive statement.  \par 
   In [DH1] and [DH2], Dunkel and Hänggi ask the question of the asymptotic behaviour of the variance, or ``mean square displacement'', of their diffusion. Indeed,  comparing to the non-relativistic case, and after numerical computations, they guess that this variance, normalised by time, should converge, to some constant for which they conjecture an empirical formula.  \par 
   
   We answer here these two questions, asked by Debbasch and Rivet in [DR], and by Dunkel and Hänggi  in [DH1], [DH2], and indeed a more general one.  We establish in fact rigourously in this article, a central limit theorem for a class of Minkowskian diffusions, to which  the two above mentioned ones, ROUP and Dunkel-Hänggi (DH) diffusion, belong.   As a consequence of our main result, we establish, for this whole class, the convergence of the normalised variance, guessed (for their particular case) in [DH1] and  [DH2]. Getting the exact expression for this limiting variance, and particularising to the (DH) diffusion, we can then invalidate and correct the wrong conjecture made in [DH1] on its expression and asymptotic behaviour (as the noise parameter goes to infinity). \par 
   To summarise the content, we begin by describing in Section \ref{sec.Minkdiff} below the class of Minkowskian diffusions we consider, which contains both ROUP and DH diffusions as particular cases.  \par 
   Then in Section \ref{sec.behav}, we present our study, leading to the following main result :   \parn 
{\bf Theorem \ref{the.main}} \quad {\it  Let  $\,(\x_t, \p_t)_{t\ge 0} = (x^i_t, p^i_t)_{1\le i\le d,\: t\ge 0}\,$ be a $\,\mathbb R^d\times \mathbb R^d$-valued diffusion solving the stochastic differential system :  \   for $\,1\le i\le d\,$, 
$$ (\star)\, \left \lbrace \begin{array}{lr}  d x_t^i  = f(|\p_t|) \, p^i_t \, d t  \\ \vspace{-2mm} \\ 
d p_t^i = -\, b(|\p_t|) \, p^i_t \, d t  + \s(|\p_t|) \Big(\beta\, [1+\eta(|\p_t|)^2]\Big)^{-1/2} [dW_t^i + \eta(|\p_t|) (p^i_t/|\p_t|)\, dw_t] & \end{array} \right. \hspace{-5mm} . $$
   Then, under some hypotheses $({\mathcal H})$ on the continuous functions $\,f$, $b\,$, $\s\,$, $\eta\,$,  the law of the process $\left( t^{-1/2} \: \x_{at} \right)_{a\geq 0}$ converges, as $\,t\to\ii\,$, to the law of $\,\left(\Sigma_\b\, \mathcal{B}_a\right)_{a \geq 0}\,$,  in $\,C(\mathbb R^+, \mathbb R^d)$ endowed with the topology of uniform convergence on compact sets of $\,\mathbb R_+\,$. Here $\,W$ and $\,\mathcal B\,$ are standard $d$-dimensional Brownian motions, $w\,$ is a standard real Brownian motion, independent of $W$, $\,\beta >0\,$ is an inverse noise or heat parameter, and $\,\Sigma_\b\,$ is a constant, displayed in Proposition \ref{pro.Martconv}  below. } \par \smallskip 
   The following is then deduced. Recall that the symbol $\,\E\,$ stands for expected or mean value, with respect to the underlying probability measure (or distribution) $\,\P$ (governing the given Brownian motions $W,w$).  \parn
{\bf Corollary \ref{cor.main}} \quad  {\it Under the same hypotheses as in the above theorem, from any starting point, the normalised variance (mean square displacement) $\,t^{-1} \: \E \left[ |\x_t|^2 \right] $ goes, as $\,t\to\ii\,$, towards $\,d\times \Sigma^2_\b\,$. }\par \smallskip 

    In Section \ref{sec.behavs}, we study the behaviours of the limiting variance $\,\Sigma^2_\b\,$, as the inverse noise parameter $\,\b\,$ goes to $0$ or to $\,\ii\,$, and we also support our result by numerical simulations. Focussing on  the particular case of the DH diffusion, for $\,d=1\,$, we get the following, which, though confirming a non-classical variance behaviour, shows up a behaviour near $0$ which differs from the one implied by the wrong guess made in [DH1] about the expression of $\,\Sigma^2_\b\,$. \smallskip \parn     
{\bf Proposition \ref{pro.estim}} \quad {\it  Consider the DH case, for $\,d=1\,$, as in [DH1]. Then, we have \par\smallskip \centerline{ ${\ds \Sigma_\beta^2 \sim 2/\b\,}$ as $\,\b\nea \infty$ ; \  and,  as $\,\b\sea 0$ :  \   ${\ds \Sigma_\beta^2 \sim {A \over\log (1/\b)} }$, \   for some explicit constant $\,A>0\,$. }   } 
\par \smallskip 

   Finally, we detail in Sections \ref{sec.ppropS} and \ref{sec.compas} two somewhat involved proofs. We thank Reinhard Sch\"afke for his kind and decisive help for the proof of  Sections  \ref{sec.psi12} and \ref{sec.psi=}.

\section{A class of Minkowskian diffusions}  \label{sec.Minkdiff}  

   Let $\mathbb R^{1,d}$, where $d \geq 1$ is an integer, denote the usual Minkowski space of special relativity. In its canonical basis, denote by $\,x = (x^{\m}) = (x^0, \: x^i) = (x^0, \:\x)$ the coordinates of the generic point, with greek indices running $\,0,..,d\,$ and latin indices running $\,1,..,d\,$. The Minkowskian pseudo-metric is given by : \  $ds^2 = |dx^0|^2 - \sum_{i=1}^d\limits |dx^i|^2$. \par
   The world line of a particle having mass $\,m\,$ is a timelike path in $\mathbb R^{1,d}$, which we can always parametrize by its arc-length, or proper time $\,s\,$. So the moves of such particle are described by a path $\,s\mapsto (x^\mu_s)$, having momentum $\,p = (p_s)$ given by : 
$$ p= (p^{\m}) = (p^0, \, p^i) = (p^0, \, \p)\, , \quad \textrm{where}\quad p^{\m}_s :=\, m\, \frac{d x^{\m}_s}{d s} \, , $$ 
and satisfying :  \qquad  ${\ds |p^{0}|^2 - \sum_{i=1}^d\limits |p^{i}|^2 = m^2\,. }$ 
\par\medskip 

   We shall consider here world lines of type $(t,\,\x(t))_{t \geq 0}\,$, and take $\,m=1$. Introducing the velocity $\,\v=(v^1, \ldots, \, v^d)$ and polar coordinates $(r,\Theta)$ by setting : 
$$ v^{i} := \frac{d x^{i}}{d t}\,,  \qquad  r := | \p | = \left( \sum_{i=1}^d\limits |p^{i}|^2 \right)^{1/2}  \quad \textrm{and } \quad \Theta := {\p\over r} =: \left( \t^1, \ldots, \, \t^d \right) \in \mathbb S^{d-1} , $$
we get at once : 
$$ p^0 = \frac{dt}{ds} = \sqrt{1+ r^2} =\left( 1- |\v|^2 \right)^{-1/2} \quad \textrm{and} \quad  p = \sqrt{1+ r^2}\, (1, \, \v) \, . $$
Thus, a full space-time trajectory 
$$ (x(t),\, p(t))_{t \geq 0}= (t, \,\x(t),\, p^0(t), \,\p(t))_{t \geq 0} $$ 
is determined by the mere knowledge of its spacial component  $(\x(t),\,\p(t))$. \par 
   We can therefore, from now on, focus on spacial trajectories $\, t\mapsto (\x_t,\,\p_t) \in \mathbb R^d\times \mathbb R^d$. 
\par\medskip 
    The Minkowskian diffusions we consider here are associated as above, to Euclidian diffusions $\, t\mapsto (\x_t, \, \p_t) = (x^i_t, p^i_t)_{1\le i\le d}\,$, which are the solution to a stochastic differential system of the following type : 
$$ (\star)\; \left \lbrace \begin{array}{lr}  d x_t^i  = f(r_t) \, p^i_t \, d t  \\ \vspace{-3mm} \\ 
d p_t^i = -\, b(r_t) \, p^i_t \, d t  + \s(r_t)\,\Big(\beta\, [1+\eta(r_t)^2]\Big)^{-1/2}\, [dW_t^i + \eta(r_t)\,\t_t^i\, dw_t] & \end{array} \right. \!\!\!\!\! , \quad \hbox{for } \; 1\le i\le d\, , $$
where $\,\W :=(W^1, \ldots , W^d)$ denotes a standard $d$-dimensional Euclidian Brownian motion, $w\,$ denotes a standard real Brownian motion, independent of $\W$, $\,\beta >0\,$ is an inverse noise or heat parameter, and the real functions $\,f$, $b\,$, $\s\,$, $\eta\,$ are continuous on $\,\mathbb R_+\,$ and satisfy the following hypotheses, for some fixed $\,\e >0$ :   
$$ ({\bf{\mathcal H}}) \qquad  \s \ge \e \;\hbox{on } \mathbb R_+ \; ; \;\,  g(r):= {2\,r\, b(r)\over \s^2(r)}\, \ge {\e} \; \hbox{ for large } r\; ;  \;\, \lim_{r\to\ii} \, e^{-\e' r} f(r) = 0\;  \hbox{ for some } \, \e'< {\ts{\beta\,\e\over 2}}\, .$$ 
\par\medskip 

   Of course, in the particular case of constant functions $f$, $b$, $\s$, and $\eta=0$, the process $(\x_t)$ is an integrated Ornstein-Uhlenbeck process. The process considered by Debbasch, Malik and Rivet ([BDR1], [BDR2], [BDR3], [DMR], [DR]), they call Relativistic Ornstein-Uhlenbeck Process (ROUP), corresponds to : 
$$ (ROUP)\qquad f(r) = b(r) = (1+r^2)^{-1/2}\, ,\; \; \s(r) = \sqrt{2}\, , \;\;\eta =0\,, \;\; g(r) =  r\, (1+r^2)^{-1/2}\,, $$
and the relativistic process considered by Dunkel and Hänggi ([DH1], [DH2]) corresponds to : 
$$ (DH)\quad f(r) = (1+r^2)^{-1/2}\, ,\; b(r) = 1\,,\; \s(r) = \sqrt{2\sqrt{1+r^2}}\, ,  \;\; \eta(r) = r\,, \;\; g(r) = r\, (1+r^2)^{-1/2}\, . $$ 

   These processes are intended to describe the motion of a point particle surrounded by a heat bath, or relativistic fluid, with respect to the rest-frame of the fluid, in which the particle diffuses.  
   The Minkowskian diffusion $(\x_t,\,\p_t)$ solving the stochastic differential system $(\star)$ is isotropic precisely when $\,\eta\equiv 0\,$ (for $d\ge 2$ ; when $d=1$, $\,\eta\,$ does not matter). If $\,\eta\not= 0\,$, the momentum $(\p_t)$ undergoes a radial drift.  
\par\medskip  

   In ([DR], Section 4), Debbasch and Rivet argue heuristically that the so-called ``hydrodynamical limit'' of their ROUP should behave in a Brownian way, and ask the question of a mathematical proof confirming such not much intuitive statement.  \par 
   In ([DH1], [DH2]), Dunkel and Hänggi ask the question of the convergence, as $\,t$ goes  to infinity, of the normalised variance (or mean square displacement) :
$$ \Sigma^2(t)\, := \, t^{-1} \: \E \left(  \sum_{i=1}^d {|x^i_t|}^2 \right) = \E\Big[ |\x_t|^2/ t\Big] . $$

   We shall answer these two questions, by means of the more general one we address, which is the asymptotic behaviour, as $\,t\to\ii\,$, of the process : 
$$ \left(\x^t_a \right)_{a \geq 0} := \, \left( t^{-1/2}\: \x_{at}\right)_{a \geq 0} \, =\, t^{-1/2} \, \left(x^1_{at}\,, \ldots, \, x^d_{at}\right)_{a \geq 0} \, , $$
where the diffusion $(\x_t,\,\p_t)_{t \geq 0}\,$ solves $(\star)$, under the hypotheses $({\cal H})$. 

\section{Asymptotic behaviour of the process ${\ds \left(\x^t_a \right)_{a \geq 0}}$} \label{sec.behav}   
\subsection{An auxiliary function $\,F\,$}   \label{sec.F} 

  Let us look for a function $\,F=(F^1,\ldots,\, F^d)\in C^2(\mathbb R^d,\mathbb R^d)\,$ such that for $\,1\le i\le d$ : 
\begin{equation} \label{eqn1} d F^i(\p_{t}) + d x^i_{t} - d  M^i_{t}=0\, , \end{equation}
for some martingale $\,\M_{t} = (M^1_t, \ldots, \, M^d_t)$, \  so that \quad 
${\ds  t^{-1/2} \: x^i_{at} = t^{-1/2} \:    M^i_{at}-  t^{-1/2} \:F^i(\p_{at}). }$ \par\medskip 

Now, Itô's Formula gives : 
$$ d F^i(\p_t) = \Bigg[ -\sum_{j=1}^d \pa_j  F^i(\p_t)\, p^j_t\, b(r_t)  + \frac{\s^2(r_t)}{2 \beta [1+\eta(r_t)^2]} \sum_{1\le j,k\le d} \Big( \delta_{jk} + \eta(r_t)^2 {\t^j_t\,\t^k_t}\Big)\,  \pa_{jk}^2F^i(\p_t) \Bigg] dt + dM^i_t \, , $$ 
with 
$$ d M_t^i= \Big(\beta\, [1+\eta(r_t)^2]\Big)^{-1/2} \s(r_t) \sum_{j=1}^d \pa_j F^i(\p_t)\, [dW_t^j + \eta(r_t)\, \t^j_t\, dw_t] . $$ 
Note that, in other words, this means that the so-called infinitesimal generator of the momentum  diffusion $(\p_t)$ is 
$$ \frac{\s^2(r)}{2\,\beta\, [1+\eta(r)^2]} \bigg( \Delta + \eta(r)^2 \sum_{1\le j,k\le d} {\t^j \t^k}\,{\partial^2\over \partial p^j \partial p^k}\bigg) - b(r) \sum_{j=1}^d\, p^j\, {\partial\over \partial p^j}\, , $$ 
$\Delta\,$ denoting the usual Euclidian Laplacian of $\mathbb R^d$. Hence a function $\,F\,$ satisfying $(1)$ must solve : 
\begin{equation} \label{eqn2}  \frac{\s^2(r)}{2\,\beta\, [1+\eta(r)^2]} \sum_{1\le j,k\le d} \Big( \delta_{jk} + \eta(r)^2\,\t^j \t^k \Big)\, \pa_{jk}^2F^i(\p) - b(r)\sum_{j=1}^d p^j\,\pa_j  F^i(\p)\,  = -\, p^i \times f(r). \end{equation}
Let us take $F^i$ of the form 
$$ F^i(\p) = \t^i \times  \psi_\beta(r) =  p^i \times \psi_\beta(r)/r \, , $$
and set for $\,r\in\mathbb R_+$ : 
\begin{equation} \label{f.ghG}  g(r) := \frac{2\, r \: b(r)}{\s^2(r)},\quad h(r) :=\frac{ 2\, r \: f(r)}{\s^2(r)} ,\quad \textrm{and} \quad G(r):=\int_0^r {g(\rho)}\, d \rho\, .  \end{equation} 
Then a direct computation shows that 
\begin{equation} \label{f.Mit}  d M_t^i = {\s(r_t)\over \sqrt{\beta [1+\eta(r_t)^2]}}  \bigg[ \psi_\beta'(r_t) d W_{t}^i + \Big[{\psi_\beta(r_t)\over r_t} -\psi_\beta'(r_t) \Big] \sum_{j=1}^d [ \d_{ij} -\t^i_{t} \t^j_{t} ] d W_{t}^j + \eta(r_t) \psi_\beta'(r_t) \t^i_t dw_t \!\bigg] ,  \end{equation} 
and that Equation $(\ref{eqn2})$ is equivalent to : 
\begin{equation} \label{f.psi}  \psi_\beta''(r) - \left( \beta\, g(r)-\frac{d-1}{r\, [1+\eta(r)^2]}\right)\psi_\beta'(r)-\frac{d-1}{r^2 [1+\eta(r)^2]}\,\psi_\beta(r) + \beta\, h(r) =0\, . \label{eqndiff3}  \end{equation}

   Note that, if $\,b \equiv f\,$,  or equivalently if $\,g \equiv h\,$, then Equation $(\ref{eqndiff3})$ admits the trivial solution $\,\psi_\beta(r)=r\,$. If $\,d=1$, Equation $(\ref{eqndiff3})$ is easily solved too.  But it is not easily solved in the general case we are considering, and not even in the case of the diffusion (DH) considered in [DH1] (and isotropically extended to higher dimensions) or in [DH2]. \par 
   However, we have the following, whose delicate proof is postponed to Section \ref{sec.ppropS}. 
\begin{prop} \label{pro.Schafke}  \  Under hypotheses $(\mathcal H)$, Equation (\ref{eqndiff3})  
admits a solution $\psi_\beta\in C^2(\mathbb R_+,\mathbb R)$  such that $\,\psi_\beta(0)=0\,$,    and \  $\,|\psi_\beta(r)| = {\mathcal O}(e^{\e' r})\,$, $\, |\psi_\beta'(r)| = {\mathcal O}(e^{\e' r})$ near infinity, for some $\,\e'< {\e\,\beta\over 2}\,$. \   Moreover, if  $\,0\le f\le b\,$, then we have $\,0\le \psi_\beta \leq Id\,$. 
\end{prop}

\subsection{Polar decomposition of the process $(\p_t)$ and equilibrium distribution $\,\nu$} \label{sec.poldec} 

   Since the diffusion $(\x_t, \: \p_t)$ solves $(\star)$, the radial process \   $ r_t = |\p_t|\,$ solves :
$$ dr_t = \left( \frac{(d-1)\,\s^2(r_t)}{2\,\beta\,[1+\eta(r_t)^2]\, r_t} - r_t\, b(r_t) \right) dt + \s(r_t)\,\beta^{-1/2}\, dB_t\,  ,  $$ 
\begin{equation} \label{f.dbB}  \hbox{with } \qquad  dB_t := (1+\eta(r_t)^2)^{-1/2} \bigg[\sum_{i=1}^d  \theta^i_t\, dW^i_t + \eta(r_t)\, dw_t\bigg] .  \end{equation} 
As $\langle B, B\rangle_t \,=\,t\,$, $\,B\,$ is a standard real Brownian motion. Consider then the angular process \   $\tilde\Theta_s = (\tilde\t^1_s, \ldots ,\, \tilde\t^d_s)\in \mathbb S^{d-1}\, $ defined by the time change \   $ \tilde\Theta_s := \Theta_{C\1(s)}$, \   i.e.  by \quad ${\ds \p_t = r_t \times \tilde\Theta_{C_t} }\,$, by means of the clock \quad ${\ds C_t =C(t) :=\int_0^t \frac{\s^2(r_s)}{\beta\, r^2_s}\, ds\, . }$ \  The process $(\tilde\Theta_s)\in\mathbb S^{d-1}$ is a spherical Brownian motion, since it solves : 
$$ d\tilde\t^i_s = \left( \frac{1-d}{2}\right) \tilde\t^i_s \, ds + \sum_{j=1}^d \left(\d_{ij}- \tilde\t^i_{s} \tilde\t^j_{s} \right) d \widetilde{W}^j_s\, ,$$
for some standard Brownian motion $\,\widetilde{W}=(\widetilde{W^1},\ldots,\widetilde{W^d}) \in \mathbb R^d$.   \  Hence the infinitesimal generator of the diffusion $(r_t, \Theta_{t}) = (r_t, \tilde\Theta_{C_t})$ is 
\begin{equation} \label{f.gener}  \mathcal A :=\, \mathcal L_r + \frac{\s^2(r)}{2\,\beta\,r^2}\: \Delta_{\mathbb S^{d-1}}\,,\quad \hbox{with} \quad  \mathcal L_r :=\frac{\s^2(r)}{2\,\beta} \,\bigg( \pa_r^2 + \bigg[ {d-1\over [1+\eta(r)^2]\, r} - \beta\, g(r) \bigg] \pa_r \bigg) . \end{equation}
Under this form, it appears that the anisotropy function $\,\eta\,$ results in a radial drift. \par  
Set for $\,r\in\mathbb R_+$ : 
\begin{equation} \label{f.mu}  \mu(r)\, := \, \exp\left[ \int_1^r {ds\over s\,[1+\eta(s)^2]}\right] \in\mathbb R_+\,  .  \end{equation} 
Note that \   $\min\{ r,1\} \le \mu(r)\le \max\{ r,1\}$, \  and that \  ${\ds 0\le r\le s\Ra 1\le {\mu(s)\over \mu(r)} \le {s\over r}\,}$. \par 

   The radial process $(r_t)$ admits the invariant measure $\,\nu(r) dr\,$, having density on 
$\mathbb R_+$ :
\begin{equation} \label{f.nu} \n(r):= \s^{-2}(r) \: \mu(r)^{d-1} \: e^{-\,\beta\, G(r)}\, . \end{equation} 

   Note that this equilibrium distribution equals the so-called J\"uttner one, in the ROUP case (we have indeed $\,G(r) = \sqrt{1+r^2}\, -1\,$ in the ROUP and DH cases). \par 
   
The hypotheses $(\mathcal H)$ ensure that $\,\n\,$ is finite, and then that the radial process $(r_t)$ is ergodic. Denoting by $\,d\Theta\,$ the uniform probability measure on the sphere $\mathbb S^{d-1}$, and setting : 
\begin{equation} \label{f.pi} \pi(dr,d \Theta) := {\textstyle \left(\int_0^\ii\nu\right)}^{-1} \times \nu(r)\, dr \, d\Theta\, , \end{equation}
it is easily seen that $\,\pi\,$ is an invariant probability measure (or equivalently : equilibrium distribution, meaning that the operator $\,\mathcal A\,$ is symmetrical with respect to $\pi$ : $\int \Phi_1\,\mathcal A \Phi_2\, d\pi = \int \Phi_2\,\mathcal A \Phi_1\, d\pi\,$ for any test-functions $\,\Phi_1,\Phi_2\,$ on $\,\mathbb R_+\times \mathbb S^{d-1}$) for the process $(r_t, \tilde\Theta_{C_t}) = (r_t, \Theta_{t})$, which is then a symmetrical ergodic diffusion on $\,\mathbb R_+\times \mathbb S^{d-1}$. 

\begin{lem} \label{lem.partF} \  For any starting point $\,\p_0= r_0 \Theta_0\,$, uniformly with respect to $\,a\ge 0\,$, we have : 
$$ t^{-1}\, \E_{\p_0} \left[  |F^i(\p_{at})|^2  \right] \longrightarrow 0\, ,  \; \hbox{ as } \; t\to\ii\,, \; \hbox{ for } \; 1\le i\le d\, . $$
\end{lem} 
\ub{Proof} \quad   Since \  $| F^i(\p)|^2 = | \t^i \times \psi_\beta(r)|^2 \leq \psi_\beta^2(r) = {\mathcal O}(e^{2\e' r})$ \  by Proposition \ref{pro.Schafke}, \   we have \par\smallskip  \centerline{ ${\ds \E_{\p_0} \left[  |F^i(\p_{at})|^2  \right] \leq C\, \E_{r_0}[ e^{2\e' r_{at}}] \,}$, for some constant $\,C\,$. } \par  
 
   Let $(Q_t)$ denote the semi-group of the radial diffusion $(r_t)$,  solution to \  
${\ds \pa_t Q_t = \mathcal L_r Q_t}\,$. It is known (see for example ([V], chapter 31)) that  $\,Q_t(r_0,r)$ is a continuous function of $(t,r)$, and that (see for example ([V], chapter 32)) $\,Q_1(r_0,r) = q_1(r)\, \nu(r),$ for some 
bounded function $\,q_1\,$. \   Hence on one hand we have : 
$$ \E_{r_0}[ e^{2\e' r_{s}} ]  \leq \sup_{0\leq s \leq 1} Q_s (e^{2\e' {\bf \cdot}})(r_0)< \infty \, , \;\hbox{ for } 0 \leq s \leq 1\, , $$
and on the other hand, by the Markov property, for $\,s\ge 1\,$ we have :
$$\E_{r_0}[ e^{2\e' r_{s}}] = Q_1 Q_{s-1} (e^{2\e' {\bf \cdot}})(r_0) = \int_{0}^\ii \E_{\rho}[ e^{2\e' r_{s-1}}]\, q_1(\rho) \nu(\rho) d\rho  \leq\,  ||q_1||_{\infty}\! \int_0^\ii e^{2\e' \rho}\,\nu(\rho)\, d\rho < + \infty\, , $$ 
by $(\mathcal H)$ and since $\,2\e'< \beta\,\e\,$. \  This shows that $\,s\mapsto \E_{r_0}[e^{2\e' r_{s}}]$ is bounded, whence the result.  $\; \diamond$ \par 

\subsection{Asymptotic study of the martingale $\,\M$} \label{sec.Mart}  
   By Formula (\ref{eqn1}) and Lemma \ref{lem.partF}, we are now left with the study of the martingale part $(\M_{t})$.  \   Recall that the coordinates $\,M^i$ of the martingale $\M$ are given by Equation (\ref{f.Mit}).

\subsubsection{Asymptotic independence of the martingales $M^i$} \label{sec.indepmart} 
\begin{lem} \label{lem.indepmart} \quad  For $\,1\le i, l\le d\,$, as $\,t\to\ii\,$ we have almost surely :  $$ \lim_{t\to\ii}\, {\langle M^i,\, M^l\rangle _t\over t} \, =\, \d_{il}\,\Sigma_\beta^2\, , \; \hbox{ with } \hskip 3mm  \Sigma_\beta^2 := {1\over \beta\,d} \bigg[ \int |\psi'_\beta|^2\, \s^2 d\pi + (d-1)\! \int \psi_\beta^2\, (1+\eta^2)\1\,Id\2\, \s^2\, d\pi\bigg] . $$ 
\end{lem} 
\ub{Proof} \quad  The computation of brackets gives easily : 
$$ \beta\, \langle M^i, \: M^l\rangle_t\, =\, \d_{il} \: S^i_t - (1-\d_{il}) \: T^{i, l}_t\, , $$ 
with 
$$\begin{array}{ll}
$$ S^i_t := $$ & \displaystyle{ \int_0^t \s^2(r_s)\,\psi_\beta'(r_s)^2\, |\t^i_{s}|^2\, ds + \int_0^t [1+\eta(r_s)^2]\1\, r_s\2 \, \s^2(r_s)\, \psi_\beta^2(r_s)\, \left(1- |\t^i_{s}|^2\right) ds }\, , \\  \\
$$ T^{i,l}_t :=$$ & \displaystyle{ \int_0^t \s^2(r_s) \Big[ \psi_\beta^2(r_s)\, [1+\eta(r_s)^2]\1\,r_s\2 -  
\psi_\beta'(r_s)^2\Big] \t^i_{s} \t^l_{s}\, ds}\, .
\end{array}$$
Setting  
$$ k^i(r,\Theta) :=  \s^2(r)\, \psi_\beta'(r)^2\, |\t^i|^2 + [1+\eta(r)^2]\1\, r\2\, \s^2(r)\,\psi_\beta^2(r)\, \big(1- |\t^i|^2\big) $$
and
$$ \ell^{i,l}(r,\Theta) := \s^2(r) \left[\psi_\beta^2(r)\, [1+\eta(r)^2]\1\, r\2 - \psi_\beta'(r)^2 \right] \t^i \t^l\, , $$ 
and noticing that these functions are $\,\pi$-integrable by Proposition \ref{pro.Schafke}, using Section \ref{sec.poldec} we can apply the ergodic theorem, to get the following almost sure convergences :  
$$ \lim_{t\to\ii} S^i_t/t\, = \int_{\mathbb R_+ \times \mathbb S^{d-1}} k^i\, d\pi\,\, ,\quad  \lim_{t\to\ii}  T^{i,l}_t  /t \, = \int_{\mathbb R_+ \times \mathbb S^{d-1}}\ell^{i,l}\, d\pi\, .$$
Now the spherical symmetry with respect to $\Theta$ implies that for $\,1 \leq i \neq l \leq d$ : 
$$ \int k^i\, d\pi\, =\, d\1 \! \int |\psi_\beta'|^2\, \s^2 d\pi + (1- d^{-1})\! \int \psi_\beta^2 \,(1+\eta^2)\1 \, Id^{-2}\,\s^2 d\pi =: \beta\,\Sigma_\beta^2\, , \quad \hbox{and} \quad  \int \ell^{i,l}\, d\pi\,= 0\, . $$ 
Hence we have got : 
$$ t^{-1}\,\langle M^i,\, M^l\rangle _t \quad \stackrel{p.s.}{\longrightarrow} \quad
\d_{il} \times \Sigma_\beta^2 = \d_{il}\, (\beta\,d)\1 \Big[ \pi\Big(\s^2 |\psi'_\beta|^2\Big) + (d-1) \,\pi\Big(\s^2 \psi_\beta^2/[(1+\eta^2) Id^2]\Big)\Big] . \;\;\diamond  $$
\par\medskip

Consider now the martingale  $\M^t$ defined by : 
$$ \M^t_{a} := (M^{1, t}_a, \ldots , \, M^{d, t}_a) :=\, t^{-1/2}\, \M_{at}\, , $$ 
and the Dambis-Dubins-Schwarz Brownian motions $\,B^{i,t}\,$, such that 
$$ M^{i, t}_a = B^{i,t}(\langle M^{i,t},M^{i,t}\rangle_a) = B^{i,t}(t\1 \langle M^{i},M^{i}\rangle_{at}) . $$  

    Applying the asymptotic Knight theorem (see for example ([RY],Theorem (2.3) and Corollary (2.4) p. 524-525)), we deduce now from Lemma \ref{lem.indepmart} the asymptotic independence of the martingales $M^i$ and $M^l$, for $1\le i \not= l\le d\,$, in the following sense. 
\begin{cor} \label{cor.indep} \  The process $(B^{1, t}, \ldots, \, B^{d, t})$ converges in law, as $\,t\,$ goes to infinity, towards a standard $d$-dimensional Brownian motion $\,\mathcal B\,$.
\end{cor} 

\subsubsection{Convergence of the finite-dimensional marginal laws} \label{convespe}

\begin{prop} \label{pro.Martconv}  \   The finite-dimensional marginal laws of the martingale $\, \M^t\,$ converge, as $\,t\,$ goes to infinity, to those of the Brownian motion $\,\Sigma_\beta \times \mathcal B$, where $\,\Sigma_\beta\,$ is the (positive) constant given by : 
\begin{equation} \label{f.Sinf1}  \Sigma_\beta^2\, = \Big[ d \int_0^\ii e^{-\,\beta\, G(r)}\,\mu(r)^{d-1} \, \s(r)\2\, dr\Big]^{-1}\times \int_0^\ii \psi_\beta(r)\, e^{-\,\beta\, G(r)}\,\mu(r)^{d-1}\, h(r)\,dr\, . \end{equation} 
Recall that $\,\psi_\beta\,$ comes from Proposition \ref{pro.Schafke}, and that we set in Formulas (\ref{f.ghG}) and (\ref{f.mu}) : 
\begin{equation} \label{f.hGbn}   G(r) = \int_0^r {g(\rho)}\, d \rho\; , \quad h(r) = \frac{ 2\, r \: f(r)}{\s^2(r)}\, ,\quad  \mu(r)\, = \, \exp\left[ \int_1^r {ds\over s\,[1+\eta(s)^2]}\right] .  \end{equation} 
\end{prop}
\ub{Proof} \quad  Fix any integer $\,N\ge 1$, positive numbers $0<a_1<\ldots <a_N\,$, and consider the vector random processes :  
$$ X^t := \Big( \langle M^{i,t}, \, M^{i,t}\rangle_{a_k}\, ,  B^{j,\,t}_s\Big)_{1\le i,j\le d\,,\, 1\le k\le N\,,\,s\ge 0}\, ,\quad  X^{\infty} := \Big( \Sigma_\beta^2 \, a_k \,,\, \mathcal B^j_{s}\Big)_{1\le i,j\le d\,,\,1\le k\le N\,,\,s\ge 0}\, . $$
By Section \ref{sec.indepmart},  $\,X^t$ converges in law, as $\,t\,$ goes to infinity, to $\,X^{\infty}$. 
   By the Skorokhod coupling theorem (see for example ([K], Theorem (4.30) p. 78)), there exist vector random processes $\,\widetilde{X}^{t}$ and $\,\widetilde{X}^{\infty}$ satisfying the identities in law :
$$ (\widetilde{X}^{t})\stackrel{d}{=} (X^{t}) \, , \quad \widetilde{X}^{\infty} \stackrel{d}{=} X^{\infty} , $$
and such that $\,\widetilde{X}^{t}$ converges almost surely to $\widetilde{X}^{\infty}$. As a consequence, we get the following convergence in distribution : 
$$ \left(  B^{i,\:t}(\langle M^{i,t},  M^{i,t}\rangle_{a_k})\right)_{1\le i\le d,\,1\le k\le N} \, 
\stackrel{d}{\longrightarrow}\, \left( \mathcal B^i(\Sigma_\beta^2 \, a_k)\right)_{1\le i\le d,\,1\le k\le N} , $$ 
or equivalently : 
$$ \left( M^{i,t}_{a_k}\right)_{1\le i\le d,\,1\le k\le N}\,  \stackrel{d}{\longrightarrow} \, 
 \Sigma_\beta \times\left( \mathcal B^i_{a_k}\right)_{1\le i\le d,\,1\le k\le N}\, .  $$
 Note that from Formulas  (\ref{f.nu}), (\ref{f.pi}), and Lemma \ref{lem.indepmart}, we get directly the following expression for $\,\Sigma_\beta$ : 
$$  \Sigma_\beta^2 =  { \int_0^\ii \psi_\beta'(r)^2\, e^{-\beta\, G(r)} \mu(r)^{d-1}\, dr + (d-1) \int_0^\ii \psi_\beta(r)^2\, e^{-\beta\, G(r)} \mu(r)^{d-1}\,[1+\eta(r)^2]\1\,r^{-2}\, dr \over  \beta d \int_0^\ii e^{-\beta\, G(r)} \mu(r)^{d-1} \s(r)\2\, dr}\,  . $$ 
    It remains to derive from this expression the expression (\ref{f.Sinf1}) of the statement for $\Sigma_\beta\,$. This is achieved as follows, integrating by parts and using Proposition \ref{pro.Schafke}, which implies that \parn  
${\ds \lim_{r\to\ii} \Big[\psi_\beta(r)\, \psi_\beta'(r)\, e^{-\beta\, G(r)} \mu(r)^{d-1} \Big] = 0\,}$, together with  Equation (\ref{f.psi}) : 
$$ \int_0^\ii \psi_\beta'(r)^2\, e^{-\beta\, G(r)}\, \mu(r)^{d-1}\, dr  $$
$$ = \Big[  \psi_\beta(r)\, \psi_\beta'(r)\, e^{-\beta\, G(r)}\, \mu(r)^{d-1}\Big]_0^\ii - \int_0^\ii  \psi_\beta(r)\, {d\over dr}\Big[\psi_\beta'(r)\, e^{-\beta\, G(r)}\, \mu(r)^{d-1} \Big] dr  $$ 
$$ =  \int_0^\ii\!\! \psi_\beta(r)  \Big[\! \Big(\! {\ts{(d-1)/r\over 1+\eta(r)^2}} - \beta g(r)\!\Big)\psi_\beta'(r)- {\ts{(d-1) \psi_\beta(r)\over r^2[1+\eta(r)^2]}} + \beta\, h(r) - \psi_\beta'(r)\Big({\ts{(d-1)/r\over 1+\eta(r)^2}} - \beta g(r)\!\Big)\! \Big] e^{-\beta G(r)} \mu(r)^{d-1} dr $$ 
$$ =\, \beta \int_0^\ii\! \psi_\beta(r)\, e^{-\,\beta\, G(r)}\,\mu(r)^{d-1}\, h(r)\,dr - (d-1)\! \int_0^\ii \!\psi_\beta(r)^2\, e^{-\beta\, G(r)}\, \mu(r)^{d-1} \,[1+\eta(r)^2]\1\, r^{-2}\, dr \, . \;\;\diamond $$

\subsubsection{Tightness} \label{sec.tension} 

\begin{prop} \label{pro.tension}  \   The family of martingales $\,\M^{t}\,$ is tight, in 
$\,C\left( \mathbb R_+, \mathbb R^d \right)$, endowed with the topology of uniform convergence on compact sets of $\,\mathbb R_+\,$. 
\end{prop}
\ub{Proof} \quad  Fix any $\,T >0$, and use the Arzel\`a-Ascoli theorem (see for example ([K], Theorem (16.5) p. 311)) : the family $\left(t^{-1/2}\, \M_{at}, \: a \in [0,T] \right)$ is tight, in $\,C\left( [0,T], \mathbb R^d \right)$, if and only if 
$$ \lim_{h\to \: 0}\, \limsup_{t \to \infty} \quad  \E \bigg[ t^{-1/2} \sup_{0 \leq a \leq  b \leq T \atop b-a \leq h} \Big|\M_{at}-  \M_{bt} \Big| \bigg] = 0\, . $$ 
Fix $\,i \in \{1, \ldots, d\}$, $h >0$, and denote by $\,n\,$ the integral part of $\, T/h\,$. There exists a standard Brownian motion $\,\widetilde{W}^i$ such that ($S^i$ being as in Section \ref{sec.Mart}) : 
$$ \frac{1}{\sqrt{ht}}\, \sup_{0 \leq a \leq  b \leq T \atop b-a \leq h}|M^i_{at}-  M^i_{bt}  |\, =\, \frac{1}{\sqrt{\beta\,ht}}\, \sup_{0 \leq a \leq  b \leq T \atop b-a \leq h}| \widetilde{W}^i \left( S^i_{bt}-S^{i}_{at} \right)|\, . $$
Setting \   $\displaystyle{ \widetilde{W^i}^*(u) = \sup_{0\leq s\leq u} |\widetilde{W}^i_s}|$, \  we have also :
$$ \frac{1}{\sqrt{ht}}\, \sup_{0 \leq a \leq  b \leq T \atop b-a \leq h} \left| \widetilde{W}^i \left(S^i_{bt}-S^{i}_{at} \right) \right|\, =\, \widetilde{W^i}^* \left( \sup_{0 \leq a \leq T-h}\, \frac{1}{ht} \left(S^i_{a+ht} -S^i_{a}\right)  \right) \leq \widetilde{W^i}^* \left(  A^i_{ht} \right) , $$
where 
$$\displaystyle{A^i_{ht} := \sup_{0 \leq j \leq n} A_{ht}^{i,j}\,, \quad \textrm{ and } \quad  A_{ht}^{i,j} := \frac{1}{ht}
\left( S^i_{(j+2)ht}- S^i_{jht}\right)}.$$
By the ergodic theorem we have (as in the proof of Lemma \ref{lem.indepmart}) the following convergence, as $\,t\to\ii\,$, valid almost surely and in $L^1$-norm as well :  \quad 
${\ds  A_{ht}^{i,j} \, {\longrightarrow} \,    2 \, \Sigma_\beta^2\, , }$ \   
which implies the uniform integrability of $\{ A_{ht}^j, \: ht \geq 1 \}$. \par  
   Otherwise, by Doob's inequality \Big(applied to the martingale ${\ds \int_0^s 1_{\{ u \leq A_{ht}\}} d \widetilde{W^i_u}\,}$\Big), we have : 
$$ \E \left[  \widetilde{W^i}^* \left( A^i_{ht} \right)  \right]\, \leq \,\Big\| \widetilde{W^i}^* \left( A^i_{ht} \right)\Big\|_2 \leq \, 2\, \sqrt{ \E [A^i_{ht}]}\, , $$
whence 
$$\E \bigg[ t^{-1/2} \sup_{0 \leq a \leq  b \leq T \atop b-a \leq h}|M^i_{at}-  M^i_{bt}  |\bigg] \leq\, 
2\,\beta^{-1/2}\, \sqrt{h \times \E [A^i_{ht}]}\, . $$
Now, as for fixed $\,h\,$ and for any $\,\la > 2 \,\Sigma_\beta^2\, $ we have :
$$ \E [A^i_{ht}]= \int_0^{\infty} \P(A^i_{ht} \geq s)ds \leq \la + \int_{\la}^{\infty} \P(A^i_{ht} \geq s)ds \leq \la+ \sum_{j=0}^n \int_{\la}^{\infty} \P(A_{ht}^{i,j} \geq s)ds\, , $$
we deduce that 
$$ \E [A_{ht}^i] \,\leq \,\la + \sum_{j=0}^n \E\left[ A_{ht}^{i,j} \times 1_{\{ A_{ht}^{i,j}\, \geq \la\}} \right] \longrightarrow \la\,, \quad \textrm{as $t$ goes to infinity.}$$
Hence, \quad ${\ds \limsup_{t \to \infty}\, \E [A_{ht}^i]\, \leq \,2 \,\Sigma_\beta^2\, , }$ \  and then 
$$ \lim_{h\to 0}\, \limsup_{t \to \infty} \quad  t^{-1/2}\, \E \bigg[ \sup_{0 \leq a \leq  b \leq T \atop b-a \leq h}|M^i_{at}-  M^i_{bt}  |\bigg] = 0\, . \;\;\diamond $$ 

\subsubsection{Main result} \label{sec.prooft} 

Gathering Formula (\ref{eqn1}), Lemma \ref{lem.partF}, and Propositions \ref{pro.Martconv} and \ref{pro.tension}, we get at once the following main result of this article.  
\begin{theo} \label{the.main}\  Let  $(\x_t, \p_t)= (x^i_t, p^i_t)_{1\le i\le d}\,$ be a $\,\mathbb R^d\times \mathbb R^d$-valued diffusion solving the stochastic differential system 
$$ (\star)\; \left \lbrace \begin{array}{lr}  d x_t^i  = f(r_t) \, p^i_t \, d t  \\ \vspace{-3mm} \\ 
d p_t^i = -\, b(r_t) \, p^i_t \, d t  + \s(r_t)\,\Big(\beta\, [1+\eta(r_t)^2]\Big)^{-1/2}\, [dW_t^i + \eta(r_t)\,\t_t^i\, dw_t] & \end{array} \right. \!\!\!\!\! , \quad \hbox{for } \; 1\le i\le d\, , $$
where $\,\W :=(W^1, \ldots , W^d)$ denotes a standard $d$-dimensional Euclidian Brownian motion, $w\,$ denotes a standard real Brownian motion, independent of $\W$, $\,\beta >0\,$ is an inverse noise or heat parameter, and the real functions $\,f$, $b\,$, $\s\,$, $\eta\,$ are continuous on $\,\mathbb R_+\,$ and satisfy the following hypotheses, for some fixed $\,\e >0$ :   
$$ ({\bf{\mathcal H}}) \qquad  \s \ge \e \;\hbox{on } \mathbb R_+ \; ; \;\,  g(r):= {2\,r\, b(r)\over \s^2(r)}\, \ge {\e} \; \hbox{ for large } r\; ;  \;\, \lim_{r\to\ii} \, e^{-\e' r} f(r) = 0\;  \hbox{ for some } \, \e'< {\ts{\beta\,\e\over 2}}\, .$$ 
  
   Then the law of the process $\left( t^{-1/2} \: \x_{at} \right)_{a\geq 0}\,$ converges, as $\,t\to\ii\,$, to the law of $\,\left(\Sigma_\beta\, \mathcal{B}_a\right)_{a \geq 0}\,$, \  
in $\,C(\mathbb R_+, \mathbb R^d)$, endowed with the topology of uniform convergence on compact sets of $\,\mathbb R_+\,$. Here $\,\mathcal B\,$ is a standard $d$-dimensional Brownian motion, and the constant $\,\Sigma_\beta\,$ is given by Formula (\ref{f.Sinf1}). 
This result holds from any starting point $(\x_0, \p_0)$ ($\p_0$ can also obey the equilibrium law $\,\pi$). 
\end{theo} 

We deduce now the result conjectured in [DH1], [DH2], and an expression of the limit. 
\begin{cor} \label{cor.main} \  Under the same hypotheses as in the above theorem, for any starting point, the normalised variance (mean square displacement) $\,t^{-1} \: \E \left[ |\x_t|^2 \right] $ goes, as $\,t\to\ii\,$, towards $\,d\times \Sigma_\beta^2\,$. 
\end{cor} 
\ub{Proof} \quad By Theorem \ref{the.main}, we have convergence in law of the random variable $\,t^{-1} |\x_t|^2$, towards $\,\Sigma_\beta^2|\mathcal{B}_1|^2$. By Formula (\ref{eqn1}) and Lemma \ref{lem.partF}, we have only to make sure that for $1\le i\le d$, the following holds : 
$$ t^{-1} \E \left[ |M^i_{t}|^2 \right] = (\beta\,t)^{-1} \: \E \left[S^i_t  \right]\longrightarrow \Sigma_\beta^2\, . $$
Now, on one hand we already noticed (recall the proof of Lemma \ref{lem.indepmart}) that, by ergodicity, we have \  ${\ds t^{-1} S^i_t \stackrel{p.s.} \longrightarrow \int k^i\, d\pi = \beta\, \Sigma_\beta^2\,}$. And on the other hand, exactly the same reasonning as in the proof of Lemma \ref{lem.partF} (to show that $\,s\mapsto \E_{r_0}[ e^{2\e' r_{at}}]$ is bounded), merely using the semi-group $(P_t)$ of the diffusion $(r_t, \Theta_{C_t})$,  solution to \  
${\ds \pa_t P_t = \mathcal A P_t}\,$, instead of the radial semi-group $(Q_t)$, 
shows that $\,s\mapsto \E_{\p_0} [k^i(\p_s)]$ is bounded. Moreover, in the same spirit, by the Markov property and by the proof of Lemma \ref{lem.indepmart}, we have : 
$$\E_{\p_0} \left[ \frac{S^i_t}{t}\right] =  \frac{1}{t} \int_0^1 P_s(k^i)(\p_0)\,  ds +
\int  \left( \frac{1}{t} \int_0^{t-1} P_s(k^i)(\p)\, ds \right) \tilde q_1(\p)\, \pi(d\p)\, , $$ 
$\,\tilde q_1\,$ being the bounded density of $\,P_1(\p_0,d\p)$ with respect to  $\,\pi(d\p)\,$.  It is clear that the first term of the right hand side goes to 0. Finally, by the  Chacon-Ornstein theorem and by dominated convergence, the second term goes indeed to $\,\beta\, \Sigma_\beta^2\,$. $\;\diamond$ 

\section{Behaviours of $\,\Sigma_\beta^2\,$, as $\,\b\sea 0\,$ and as $\,\b\nea \ii$} \label{sec.behavs} 

   Theorem \ref{the.main} and Corollary \ref{cor.main} show up the interest of the limiting constant $\,d\times\Sigma_\beta^2\,$.   Recall then from Sections \ref{sec.Minkdiff} and \ref{sec.F}  that the processes considered by ([BDR1], [BDR2], [BDR3], [DMR], [DR]) and by ([DH1], [DH2]), correspond respectively to : 
$$ (ROUP)\quad h(r) = g(r) =  r\, (1+r^2)^{-1/2}\, ,\; \; G(r) = \sqrt{1+r^2}-1\, ,\;\;\eta =0\,, \; \; \s(r)^2 = 2\, ,\; \;\psi_\beta(r) = r\,\, , $$
$$ (DH)\;  h(r) = {r\over 1+r^2}\, ,\, g(r) = \mu(r) = {r\over \sqrt{1+r^2}}\, ,\, G(r) = \sqrt{1+r^2}-1\,,\, \s(r)^2 = 2\,\sqrt{1+r^2} , \, \eta(r) = r\, , $$ 
\parn 
for some positive (noise or heat) inverse parameter $\,\b\,$. It is natural to wonder, as in [DH1], how behaves the limiting variance $\,\Sigma_\beta^2\,$, as $\,\b\sea 0\,$ and as $\,\b\nea \ii\,$. \par 
 
    In the ROUP case, we have simply \  $d\times\Sigma_\beta^2 = 2d/\b\,$. \  The variance behaviour is Euclidian. \par 
    In the DH case of [DH1], [DH2], we have by Formula (\ref{f.Sinf1}) : 
$$ d\times  \Sigma_\beta^2\, =\, 2\times {{\ds \int_0^\ii \psi_\beta(r)\, e^{-\beta \sqrt{1+r^2}}\, (1+r^2)^{-(d+1)/2}\, r^{d}\, dr }\over  {\ds \int_0^\ii e^{-\beta \sqrt{1+r^2}}\, (1+r^2)^{-d/2}\, r^{d-1}\, dr} } \, . $$ 
 Note that the precise value of $\,\psi_\beta\,$ is given in Section \ref{sec.psi=} :   
$$\psi_\beta(r) = \zeta_1(r) \int_0^r \zeta_2(\rho) w_\beta(\rho)\1 h(\rho) d \rho + \zeta_2(r) \int_r^{\infty} \zeta_1(\rho) w_\beta(\rho)\1 h(\rho) d \rho\, , $$ 
with functions $\,\zeta_1, \zeta_2, w_\beta\,$ given in Section \ref{sec.psi12}. \par

   In [DH1], for $\,d=1\,$,  after numerical simulations, Dunkel and Hänggi conjecture that   $\,\Sigma_\beta^2\,$ could be equal to \   ${ \frac{2}{2+\b}}\,$. The expression we got above for $\,\Sigma_\beta^2\,$ invalidates this conjecture, and, even the asymptotic behaviour near 0 it implies. However, it is true that a non-classical variance behaviour occurs. We have indeed the following, whose technical proof is postponed to Section \ref{sec.compas}. 
\begin{prop} \label{pro.estim} \   Consider the DH case, for $\,d=1\,$, as in [DH1]. Then, we have \par\smallskip \centerline{ ${\ds \Sigma_\beta^2 \sim 2/\b\,}$ as $\,\b\nea \infty$ ; \  and,  as $\,\b\sea 0$ :  \   ${\ds \Sigma_\beta^2 \sim {A \over\log (1/\b)} }$, \   for some explicit constant $\,A>0\,$. }   
\end{prop}

\subsection{Numerical Simulations} \label{sec.numsim} 

   To confirm the validity of our estimates in Proposition \ref{pro.estim}, invalidating the conjecture of [DH1], we performed numerical simulations relating to the DH diffusion, in the case $\,d=1\,$. We used the Monte-Carlo method,  with $\,N=1000\,$ simulations. For different values of $\,\b$ (from $10^{-5}$ to $10^6$), we computed  $\,\x_j(t)_{ j=1..N}\,$ for $\,0\le t\le T=1000$, and then the quantity :
$$\overline{\x^2}(T)=\frac{1}{N} \sum_{j=1}^{N} \x_j^2(T) . $$

  The following diagram represents our results in logarithmic coordinates.  Thus, the horizontal axis represents $\,\log(1/\b)$, the points \textcolor{blue}{$*$} represent the simulated values 
$\,\log(\overline{\x^2}(T)/T)$ in function of $\,\log(1/\b)$.

 The straight line corresponds to the Euclidian behaviour, the continuous curve to the function $\,\b \to 2/ (\b+2)$, and the dashed curve corresponds to 
a decrease in $\,\log(1/\b)^{-1}$ for small $\,\b\,$.

\par
\vspace{0.5cm}
\hspace{-0.2cm}
\includegraphics{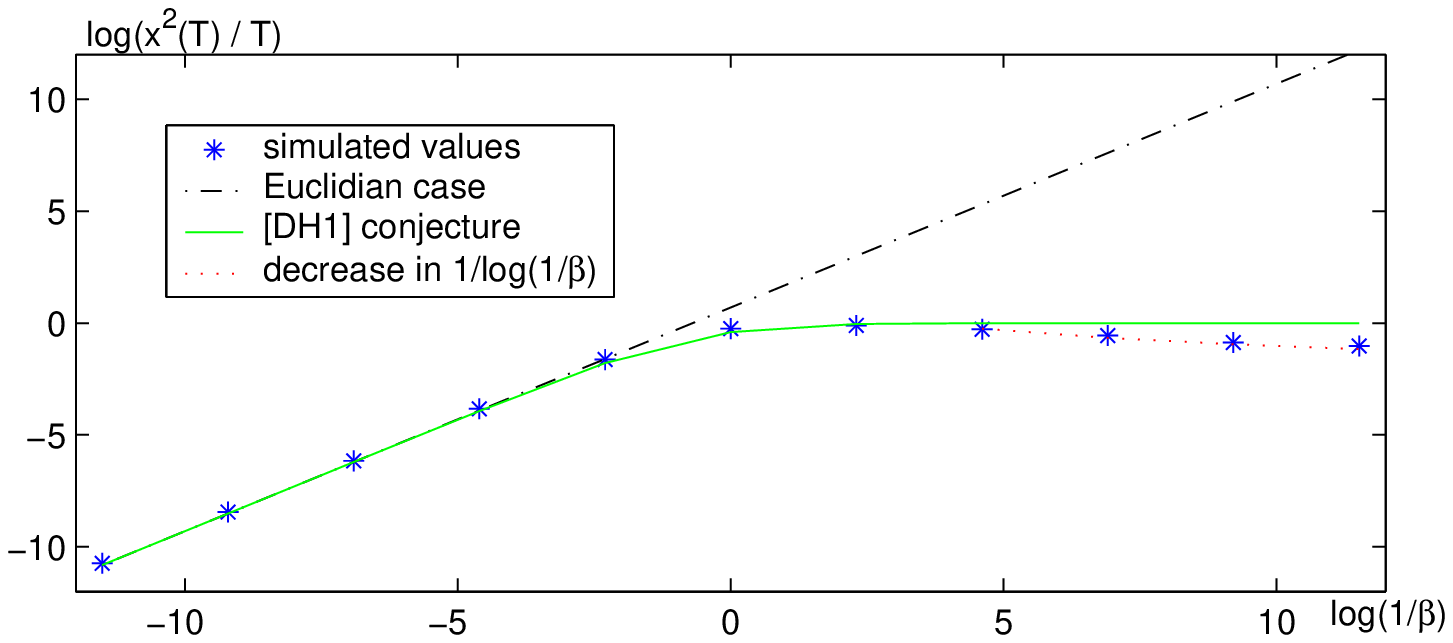} \par 
\vspace{0.5cm}

   These simulations confirm the Euclidian behaviour of the DH diffusion as $\,\b >>1\,$. For small $\,\b\,$, the expression conjectured in [DH1] is a good approximation as long as 
 $\,\b > 1/10$ ; however, for smaller $\,\b\,$, a divergence appears clearly. On the contrary, the $\log(1/\b)^{-1}$-like asymptotic behaviour of the limit, which we established above, appears as  confirmed. 

\subsubsection{The program used for the simulations} \label{sec.program} 

\begin{minipage}{18cm} \par
\noindent
\texttt{function res=asymp(N,h,D,T)}  \   (written in ``matlab'') \\
\noindent \texttt{\textcolor{blue}{\it $N$ is the iteration number in the Monte
Carlo method. The process $x(t)$ is simulted on $[0, T]$, \\ 
with mesh $\,h$. Different values for $\b$ have been tested. }}\\
\noindent \texttt{\textcolor{red}{\it Initialisation. Arrays $\,p\,$ and $\,x\,$ contain the values of $\,p(t)$ and $\,x(t)$ for $0\le t \le T$}} \\
     \texttt{t=0:h:T;}\  \texttt{n=length(t);} \  \texttt{r=[];} \\ 
\noindent
\texttt{for k=1:1:N} \   \texttt{p=zeros(1,n);} \     \texttt{x=zeros(1,n);} \\ 
\noindent 
\texttt{\textcolor{red}{\it Simulation of Brownian motion}} \\
     \texttt{u=randn(1,n);} \    \texttt{W=sqrt(2*D*h)*u;} \\
\noindent
\texttt{\textcolor{red}{\it Simulation of processes $\,p(t)$, and $\,x(t)$ by integration}} \\
     \texttt{for j=1:1:n-1} \par
          \texttt{gam=sqrt(1+p(1,j).*p(1,j));} \par
          \texttt{p(1,j+1)=p(1,j)-(p(1,j))*h + sqrt(gam)*W(1,j);}\par
          \texttt{x(1,j+1)=x(1,j)+(p(1,j)/gam)*h;}\par
\noindent   \texttt{end}\\
\texttt{\textcolor{red}{\it The $N$ simulations of $\,x(t)$ are placed in the array $\,r$}} \\
     \texttt{r=[r ; x];} \   \texttt{end} \\
\noindent
\texttt{\textcolor{red}{\it Computation of the mean of $\,\texttt{x}^{2}(\texttt{T})$, normalised by $\,T$}} \\
     \texttt{car=r.*r;} \   \texttt{limite=mean(car);}\   \texttt{res=limite(n)/T;} \\ 
\texttt{\textcolor{red}{\it end of program.}}
\end{minipage}

\section{Proof of Proposition \ref{pro.Schafke}} \label{sec.ppropS} 

   We are indebted to Reinhard Sch\"afke for this proof, who kindly indicated to us how to proceed for Sections \ref{sec.psi12} and \ref{sec.psi=} below. We thank him warmly. \   Consider first the homogeneous equation associated to $(\ref{eqndiff3})$ :
\begin{equation} \zeta''(r) + \left(\frac{d-1}{r\, [1+\eta(r)^2]} - \beta\, g(r)\right)\zeta'(r)-\frac{d-1}{r^2 [1+\eta(r)^2]}\,\zeta(r) = 0\, . \label{eqndiff3homo} \end{equation}
It has a pole of order 2 at 0 (except for $d=1$), and a pole at infinity. Using the fixed point method, we construct two solutions $\,\zeta_1\,$ and $\,\zeta_2\,$ of Equation $(\ref{eqndiff3homo})$, bounded respectively near infinity and near 0. Using these two solutions of the homogeneous equation, 
a solution $\,\psi_\beta\,$ to $(\ref{eqndiff3})$ is then deduced, which vanishes at 0. Finally, we establish the wanted control on $\,\psi_\beta , \psi_\beta'\,$. \par 
   Recall from Formula (\ref{f.mu}) that we set : \quad 
${\ds   \mu(r)  = \exp\left[ \int_1^r {ds\over s\,[1+\eta(s)^2]}\right] , }$ \   
so that $\,\mu\,$ increases and \   $\min\{ r,1\} \le \mu(r)\le \max\{ r,1\}$, \  and \  $0\le r\le s\Ra 1\le {\mu(s)\over \mu(r)} \le {s\over r}\,$. 

\subsection{Constructing solutions to the homogeneous equation $(\ref{eqndiff3homo})$} \label{sec.psi12} 

\subsubsection{Constructing a solution $\,\zeta_1\,$ to $(\ref{eqndiff3homo})$, bounded near $\infty$} \label{sec.zeta1} 

Using hypotheses $(\mathcal H)$, fix $\,\varepsilon > 0$ and $r_0\ge 1\,$ such that $\,g \geq \varepsilon\,$ on $[r_0, \infty[\,$.
For $\,r \geq r_0\,$, set 
$$ \la(r) := \int_r^{\infty} \mu(\rho)^{1-d}\,e^{{\beta}\,G(\rho)} \left[ \int_{\rho}^{\infty} e^{-\beta\,G(s)} \,\mu(s)^{d-1}\,\, [1+\eta(s)^2]\1\, s\2\, ds \right] d \rho \, . $$
We have 
$$ \la(r)\, \leq\, \int_r^{\infty} \left[ \int_{\rho}^{\infty} e^{- \beta\,\varepsilon (s- \rho)}\, \Big[{s\over \rho}\Big]^{d-1}\,s\2\, ds \right] d \rho\, = \int_r^{\infty} \left[ \int_{0}^{\infty} e^{- \beta\,\varepsilon\, s}\, (1+s/\rho)^{d-3}\, ds \right] \rho\2\, {d \rho} $$
$$ \le\,  {1\over r} \int_0^{\infty}  e^{- \beta\,\varepsilon\, s}\, \max\{ 1,(1+s/r_0)^{d-3}\}\, ds \,  = {\mathcal O}(1/r)  . $$
As  $\,r\to\ii\,$, $\la(r)$ decreases to 0, so that (up to increase $\,r_0$) we can suppose that $\,\la(r_0) \leq 1/(2d)$.
On $[r_0, \infty[\,$, let us define by induction on $\,n \in \mathbb N\,$ the functions : $\,\varphi_0 \equiv 1$, and 
$$ \varphi_{n+1}(r):= 1+ (d-1) \int_r^{\infty} \mu(\rho)^{1-d}\, e^{\beta\,G(\rho)} \left[ \int_{\rho}^{\infty}  e^{-\beta\,G(s)}\,\mu(s)^{d-1}\,\varphi_n(s)\, \, [1+\eta(s)^2]\1\, s\2\, ds \right] d\rho\,. $$
We have for $\,r \geq r_0$ :
$$ 1\leq \varphi_{n+1}(r)\, \leq\, 1 + (d-1)\, \|\varphi_n\|_{L^{\infty}[r_0, \infty[} \times \la(r) \le 1+\5 \,   \|\varphi_n\|_{L^{\infty}[r_0, \infty[} \, , $$
whence \quad 
${\ds 1 \leq\, \varphi_{n}\, \leq\, \|\varphi_n\|_{L^{\infty}[r_0, \infty[}  < 2 \, , \; \textrm{for any $\,n \in \mathbb N$}\, . }$  \   
Then similarly : 
$$ \|\varphi_{n+1} - \varphi_n\|_{L^{\infty}[r_0, \infty[} \,\leq\, \|\varphi_{n} - \varphi_{n-1}\|_{L^{\infty}[r_0, \infty[} \times (d-1) \la(r_0)\, , $$
which allows to apply the fixed point method, to get \  $\displaystyle{\zeta_1 := \lim_{n \to\infty \atop L^{\infty}[r_0, \infty[} \varphi_n}$, \  which satisfies  
\begin{equation}   1\le \zeta_1(r) = 1 + (d-1)\! \int_r^{\infty} \! \mu(\rho)^{1-d}\, e^{\beta\, G(\rho)}\! \left[ \int_{\rho}^{\infty}\! e^{-\beta\, G(s)}\,\mu(s)^{d-1} \zeta_1(s) [1+\eta(s)^2]\1 s\2 ds \right]\! d\rho\, . \label{ptfixe1}\end{equation} 
In particular, as $\,r\to\ii\,$ we have : 
$$ \zeta_1(r) \leq  1 + (d-1)\, \|\zeta_1\|_{L^{\infty}[r_0, \infty[}\, \la(r) \longrightarrow 1\, , $$
hence \   $\lim_{r \to \infty}\limits  \zeta_1(r)=1\,$,  and 
$$ \zeta_1'(r) = (1-d)\,  \mu(r)^{1-d}\, e^{\beta\,G(r)}\int_{r}^{\infty} e^{-\beta\, G(s)}\,\mu(s)^{d-1}\, \zeta_1(s)\, [1+\eta(s)^2]\1\, s\2\, ds\, < 0\, , $$
then 
$$ \zeta_1''(r) + \left(\frac{d-1}{r\, [1+\eta(r)^2]}- \beta\,g(r)\right)\zeta_1'(r)-\frac{d-1}{r^2  [1+\eta(r)^2]}\,\zeta_1(r) = 0\, . $$
This solution can be continued over the whole $\,\mathbb R_+^*\,$, yielding $\,\zeta_1\,$ still satisfying $(\ref{eqndiff3homo})$ and $(\ref{ptfixe1})$ on $\,\mathbb R_+^*\,$. We have also  
$\,\lim_{0}\limits\, \zeta_1 = +\infty\,$. 

\subsubsection{Constructing a solution $\,\zeta_2\,$ to $(\ref{eqndiff3homo})$, bounded near $0$}
\label{sec.zeta2} 

   For $\,r \in [0,1]$,  set : \  
$\,{\ds \La(r):= \beta \int_0^r \mu(\rho)^{-d-1}\, \rho\2\, e^{\beta\,G(\rho)}\! \left[\int_0^{\rho} e^{-\beta\,G(s)}\, \mu(s)^{d+1}\, |g(s)|\, s\,{ds} \right]\! d \rho\,  . }$ \parn   
We have \quad ${\ds 0\le \, \La'(r) \le\,  \beta\, e^{2\beta \int_0^r|g|}\! \int_0^1 |g(rs)|\,s\,ds\, \longrightarrow 0 \quad \textrm{as $\,r \to 0\,$}}$, by  hypotheses $(\mathcal H)$. We can then 
fix $\,r_1 \in\, ]0,1]$ such that $\,\La(r_1) \leq 1/2\,$. On $]0, r_1]$,  let us define by induction on $\,n \in \mathbb N\,$ the functions : $\,\phi_0 \equiv 1$, and 
$$ \phi_{n+1}(r) := 1 + \beta \int_0^r \mu(\rho)^{-d-1}\, \rho\2\, e^{\beta\,G(\rho)}\, \left[\int_0^{\rho} e^{-\beta\,G(s)} \mu(s)^{d+1}\, g(s)\,\phi_n(s)\, s\, ds \right] d \rho\, . $$
We have \   $\phi_n \in C^2(]0, r_1])$, \  $\|\phi_{n+1}\|_{L^{\infty}]0, r_1]} \le  1 + \La(r_1) \|\phi_{n}\|_{L^{\infty}]0, r_1]}\,$ so that \  $\|\phi_{n}\|_{L^{\infty}]0, r_1]} < 2\,$,  and 
$$ \|\phi_{n+1}- \phi_{n}\|_{L^{\infty}]0, r_1]}\, \leq\, \|\phi_{n}- \phi_{n-1}\|_{L^{\infty}]0, r_1]}\times  \La(r_1) , $$
which allows to apply the fixed point method, to get \quad   
$\displaystyle{\widetilde{\phi} :=\lim_{n \to \infty \atop  L^{\infty}]0, r_1]} \phi_n}\,$, \  which satisfies :  
\begin{equation} \widetilde{\phi}(r) = 1+ \beta\! \int_0^r\! \mu(\rho)^{-d-1}\, \rho\2\, e^{\beta\, G(\rho)}\! \left[\int_0^{\rho}\! e^{-\beta\, G(s)}\, \mu(s)^{d+1}\, g(s)\, \widetilde{\phi}(s)\, s\, ds \right]\! d \rho\, = 1+ {\mathcal O}[\Lambda(r)] \textrm{} \label{ptfixe2}\, \end{equation}
for any $\,r \in\, ]0, r_1]$. \   Hence, 
$$  \widetilde{\phi}' (r) = \beta\,\mu(r)^{-d-1}\, r\2\, e^{\beta\,G(r)}\! \int_0^{r}\! e^{-\beta\,G(s)}  \mu(s)^{d+1}\, g(s)\, \widetilde{\phi}(s)\, s\, ds  = {\mathcal O}[\Lambda'(r)]
\longrightarrow 0\quad  \textrm{ as } r \to 0\, . $$ 
Therefore,  $\,\widetilde{\phi}(0)=1$, $\,\widetilde{\phi}'(0)=0\,$, \  and for any $\,r \in\, ]0, r_1]$ :
\begin{equation}  \widetilde{\phi}''(r) + \left[ \frac{d+1}{r\, [1+\eta(r)^2]} + \frac{2}{r} - \beta\, g(r)  \right]\widetilde{\phi}'(r) - \beta\, \frac{g(r)}{r}\, \widetilde{\phi}(r)=0\, .  \label{ptfixe21} \end{equation}
This function $\,\widetilde{\phi}\,$ can be continued on the whole $\,\mathbb R_+\,$, into a function $\,\widetilde{\phi}\,$ satisfying still Equations $(\ref{ptfixe2})$ and $(\ref{ptfixe21})$.  \  
   Set now \quad $\zeta_2(r) := r\, \widetilde{\phi}(r)$. \   It is immediate that $\,\zeta_2\,$ 
solves $(\ref{eqndiff3homo})$ on $\,\mathbb R_+\,$, and satisfies :  \parn 
\centerline{$\,\zeta_2(0)=0\,$, $\,\zeta_2'(0)=1\,$.} 

\subsubsection{The Wronskian $\,w_\beta\,$ of $\,\zeta_1, \zeta_2$} \label{sec.wronsk} 

Consider the  Wronskian :
$$ w_\beta := \zeta_1 \zeta_2'-\zeta_1' \zeta_2 \quad \textrm{on $\,\mathbb R_+^*$\,}. $$
We have 
$$ w_\beta'= \zeta_1 \zeta_2''- \zeta_1'' \zeta_2 = \left( \beta\, g- \frac{d-1}{[1+\eta^2]\,Id} \right) \times w_\beta\, ,$$
so that \   
$$ w_\beta(r) = a_\beta\, \mu(r)^{1-d}\, e^{\beta\,G(r)}\, \quad \hbox{for any $\,r >0$} $$ 
and for some constant $\,a_\beta\,$. As $\,\zeta_1 \geq 1\,$, $\,\zeta_2' >0\,$, $\,\zeta_1' < 0\,$, $\,\zeta_2 >0\,$  near $\,0\,$,
we must have  $\,a_\beta >0\,$.

\subsection{Constructing a solution $\,\psi_\beta\,$ to Equation $(\ref{eqndiff3})$ on $\,\mathbb R_+$}  \label{sec.psi=} 

   For any continuous function  $\,k\,$ on $\,\mathbb R_+\,$, such that \  
${\ds  \lim_{r\to\ii} \, e^{-\e' r} k(r) = 0\,}$  for some $\,\e'< {\ts{\beta\,\e\over 2}}\,$, \ and for any $\,0 <r < \infty\,$,  \  
set :    
$$ \Psi(k)(r) := \zeta_1(r) \int_0^r \zeta_2(\rho)\, w_\beta(\rho)\1\, k(\rho)\, d \rho + \zeta_2(r) \int_r^{\infty} \zeta_1(\rho)\, w_\beta(\rho)\1\, k(\rho)\, d \rho\,  . $$ 
Note that $\,\Psi(k)$ is well defined ; we have indeed, using that $\,g\ge\e\,$ on $[r_0,\ii[$ : 
$$ \int_{r_0}^{\infty} \,\zeta_1(\rho)\,w_\beta(\rho)\1 |k(\rho)|\, d \rho\, =\,
{\mathcal O}(1) \int_{r_0}^{\infty} \mu(\rho)^{1-d}\, e^{-\beta\, G(\rho)}\, e^{\beta\, \e\,\rho/2}\, d \rho\, <  \infty\, . $$
Note also that by $(\mathcal H)$, we can take in particular $\,k=\beta\,h\,$ (recall Formula (\ref{f.ghG}) defining $\,h$). \parn 
Moreover, again for $\,0 <r < \infty\,$, we have : 
$$ \Psi(k)'(r) = \zeta_1'(r) \int_0^r \zeta_2(\rho)\, w_\beta(\rho)\1 k(\rho)\, d \rho + \zeta_2'(r) \int_r^{\infty} \zeta_1(\rho)\, w_\beta(\rho)\1 k(\rho)\, d \rho\, , $$
and 
$$ \Psi(k)''(r) = \zeta_1''(r) \int_0^r \psi_2(\rho)\, w_\beta(\rho)\1 k(\rho)\, d \rho + \zeta_2''(r) \int_r^{\infty} \zeta_1(\rho)\, w_\beta(\rho)\1 k(\rho)\, d\rho - k(r)\, , $$
so that $\,\Psi(k)\,$ solves on $\,\mathbb R_+^*\,$ Equation $(\ref{eqndiff3})$, with $\,k\,$ instead of  $\,\beta\,h\,$. Near 0, we have \  $\zeta_2(r) \sim r\,$. Otherwise, noticing that \  $(\zeta_1/\zeta_2)' = - w_\beta/\zeta_2^2\,$,  we have for $\,r>0$ :
$$ \int_r^1 \frac{w_\beta}{\zeta_2^2} = \frac{\zeta_1}{ \zeta_2}(r) - \frac{\zeta_1}{\zeta_2}(1), \quad i.e., \quad \zeta_1(r) = \frac{\zeta_1}{\zeta_2}(1)\,\zeta_2(r) + \zeta_2(r) \int_r^1 \frac{w_\beta}{\zeta_2^2}\, .  $$
Hence, near 0 we have :  $\,{\ds\zeta_1(r)\sim r \, a_\beta \int_r^1 \mu(s)^{1-d} s^{-2} ds \le  {a_\beta}\, \mu(r)^{1-d}}$, and then  \  $\Psi(k)(r) = {\mathcal O}(r)$. \  In particular, we have $\,\Psi(k)(0)=0\,$. \  
Using $\,\zeta_1(s) = {\mathcal O}\Big(\mu(s)^{1-d}\Big)$ in the expression of $\,\zeta_1'$ (recall Section \ref{sec.zeta1}), we get at once \   $|\zeta_1'(r)| = {\mathcal O}\Big(\mu(r)^{1-d} /r\Big)$ near $0\,$.  
\      We have therefore near 0 : 
$$ |\zeta_1'(r)|  \int_0^r \zeta_2(\rho)\, w_\beta(\rho)\1 k(\rho)\, d \rho\, =\,  {\mathcal O}\Big(\mu(r)^{1-d} /r\Big) \int_0^r \mu(\rho)^{d-1}\, \rho\,d\rho\,  =\, {\mathcal O}(r) , $$ 
and 
$$ \zeta_2'(r) \int_r^{\infty} \zeta_1(\rho)\, w_\beta(\rho)\1 k(\rho)\, d \rho\, =\, {\mathcal O}(1),$$
whence \quad   ${\ds \Psi(k)'(0) = \int_0^{\infty} \zeta_1(\rho)\, w_\beta(\rho)\1 k(\rho)\, d \rho\, \in \mathbb R}\,$. \   Using $(\ref{eqndiff3})$, we get also $\,\psi_\beta''(0) \in \mathbb R\,$. \par 
\smallskip 
     Setting \   $\psi_\beta := \beta\,\Psi(h)$, we have thus \  $\,\psi_\beta\in C^2(\mathbb R_+)$,  $\psi_\beta(0)=0\,$, and $\,\psi_\beta\,$ solves $(\ref{eqndiff3})$ on $\,\mathbb R_+\,$. 

\subsection{Estimates for $\,\psi_\beta$ and $\,\psi_\beta'$ near $\,\ii$} \label{sec.integp} 

  Recall from Section \ref{sec.psi=} that 
$$ \psi_\beta(r) = \beta\, \zeta_1(r) \int_0^r \zeta_2(\rho)\, w_\beta(\rho)\1 h(\rho)\, d \rho + \beta\, \zeta_2(r) \int_r^{\infty} \zeta_1(\rho)\, w_\beta(\rho)\1 h(\rho)\, d \rho $$ 
and 
$$ \psi_\beta'(r) = \beta\, \zeta_1'(r) \int_0^r \zeta_2(\rho)\, w_\beta(\rho)\1 h(\rho)\, d \rho + \beta\, \zeta_2'(r) \int_r^{\infty} \zeta_1(\rho)\, w_\beta(\rho)\1 h(\rho)\, d \rho\, . $$
Near infinity, we have on one hand : \ $\zeta_1\sim 1\,$ and  (by Section \ref{sec.zeta1}) 
$${\ds |\zeta_1'(r)| = {\mathcal O}\Big(\mu(r)^{1-d}\Big)\!\! \int_r^\ii\!  e^{-\beta\,\e\,(s-r)}\, \mu(s)^{d-1}\, s^{-2}\, ds = {\mathcal O}(r^{-2})\! \int_0^\ii  e^{-\beta\,\e\,s} \, (1+s/r)^{d-1}\, ds = {\mathcal O}(r^{-2})} . $$   
On the other hand, by $(\mathcal H)$ we have also near infinity : 
$$ w_\beta(r)-\frac{2}{\beta\,\varepsilon}\, w_\beta'(r)= \left(1-\frac{2}{\varepsilon}\, g(r) + \frac{2\,(d-1)}{\beta\,\varepsilon\,[1+\eta(r)^2]\, r} \right) w_\beta(r) < 0\, , $$
so that \qquad ${\ds\int_{r_0}^r w_\beta\, \leq \frac{2}{\beta\,\varepsilon}\, \left( w_\beta(r)-w_\beta(r_0) \right)} = {\mathcal O}(w_\beta(r))$. \parn 
Noticing that \  $(\zeta_2/\zeta_1)' = w_\beta/\zeta_1^2\,$, we have for any $\,r>0$ : \quad 
${\ds   \zeta_2(r) = \frac{\zeta_2}{\zeta_1}(1)\,\zeta_1(r) + \zeta_1(r) \int^r_1 \frac{w_\beta}{\zeta_1^2}\, . }$ \parn 
This implies that \quad $ \zeta_2(r) =  {\mathcal O}(w_\beta(r)) $ \   near infinity. \   
Therefore, by $(\mathcal H)$ there exists an $\e'< \e\beta/2\,$ such that :
$$ \int_0^r \zeta_2(\rho)\, w_\beta(\rho)\1 h(\rho)\, d \rho =  \int_0^r {\mathcal O}(e^{\e' \rho})\, d \rho = {\mathcal O}(e^{\e' r})\, . $$ 

   To control the second integral in the expressions for  $\,\psi_\beta$ and $\,\psi_\beta'\,$, observe similarly that :  
$$ w_\beta(r) \int_r^{\infty} \zeta_1(\rho)\, w_\beta(\rho)\1 h(\rho)\, d \rho\, =\,  \mu(r)^{1-d} \int_r^{\infty}  e^{-\beta [G(\rho)-G(r)]}\, \mu(\rho)^{d-1}\, {\mathcal O}(e^{\e' \rho})\, d \rho $$
$$ =\,  {\mathcal O}(e^{\e' r})\, \mu(r)^{1-d}\!  \int_r^\ii e^{- \e\,\beta (\rho-r)/2}\,\mu(\rho)^{d-1}\, d \rho\, =\,  {\mathcal O}(e^{\e' r})\! \int_0^\ii e^{- \e\,\beta\,\rho /2}\, (1+\rho /r)^{d-1}\, d \rho\, = \, {\mathcal O}(e^{\e' r})\, . $$ 
   Then by definition of $\,w_\beta\,$, we have \quad ${\ds \zeta_2'= \frac{w_\beta}{\zeta_1} + \frac{\zeta_2}{\zeta_1}\,\zeta_1'}\,$, \quad   whence \quad  $\zeta_2' \sim w_\beta$ \   near infinity.   \par\smallskip 
   
   As a conclusion, gathering the above we have indeed, for some $\,\e'< \e\beta /2$ : \  $\, |\psi_\beta(r)| = {\mathcal O}(e^{\e' r})\,$ and $\, |\psi_\beta'(r)| = {\mathcal O}(e^{\e' r})\,$, \  for large $\,r\,$. \par 

\subsection{We have \  $\,\beta\, \Psi(g) = Id\,$ on $\,\mathbb R_+\,$, \  and  \   $\,0\le \psi_\beta\le Id\,$ if $\,0\le f\le b$} \label{sec.encpsi} 

   Set $\,\tilde{\psi}_\beta := \beta\, \Psi(g)$, and note that the identical function $\,Id\,$ solves $(\ref{eqndiff3})$ if $\,g=h\,$. Hence, the function $\,r \mapsto \widetilde{\psi}_\beta(r)-r\,$ solves the homogeneous equation $(\ref{eqndiff3homo})$, so that 
$$ \widetilde{\psi}_\beta(r)-r = c\: \zeta_2(r) + c'\, \zeta_1(r)\, , \quad \hbox{ for some real constants } \; c,c'\, \; \hbox{ and for any }\,  r>0\, . $$
As $\,\widetilde{\psi}_\beta(0)=0\,$, we must have $\,c'=0\,$. Then, as \   
${\ds \zeta_2(r)\,\sim \int_1^r w_\beta\, \gg\, e^{\beta\, G(r)/2} \gg r \, }$ for large $\,r\,$,  
we must have  $\,c \geq 0\,$. Otherwise, integrating by parts, near infinity we have : 
$$ \beta\! \int_r^{\infty}\! \mu(\rho)^{d-1} e^{-\beta\, G(\rho)}  g(\rho) d \rho = \mu(r)^{d-1} e^{-\beta\, G(r)} + {\mathcal O}(r\1)\!\! \int_r^{\infty}\! e^{-\beta\, G(\rho)} \mu(\rho)^{d-1}\, d \rho\,  \sim\, \mu(r)^{d-1}\, e^{-\beta\, G(r)}, $$
so that 
$$ \beta \int_r^{\infty} \zeta_1(\rho)\, \mu(\rho)^{d-1}\, e^{-\beta\, G(\rho)}  g(\rho)\, d \rho\, \sim \beta\int_r^{\infty} \mu(\rho)^{d-1}\, e^{-\beta\, G(\rho)}  g(\rho)\, d \rho\, \sim \mu(r)^{d-1}\, e^{-\beta\, G(r)} . $$
Hence, by definition of $\,\tilde{\psi}_\beta\,$ and $\,\Psi$ and by the above, we have near infinity :
$$ \tilde{\psi}_\beta(r) = {\mathcal O}(1) \int_0^r g(\rho)\, d \rho + {\mathcal O}(w_\beta(r))\, \mu(r)^{d-1}\, e^{-\beta\, G(r)} = {\mathcal O}(G(r)+1) = o(e^{\beta\, G(r)/2}) , $$ 
whence 
$$ c \, e^{\beta\, G(r)/2} = o\Big(c \, \zeta_2(r)\Big) = o\Big(\widetilde{\psi}_\beta(r)-r\Big) = o\Big( e^{\beta\, G(r)/2}\Big), \quad \mbox{ which forces} \quad c=0\, . $$ 
Therefore $\,\tilde{\psi}_\beta(r)=r\,$ on $\,\mathbb R_+\,$, as wanted. \  Finally, if $\,0\le f\le b\,$, then $\,0\le h\le g\,$, so that \  $\Psi(h)\ge 0\,$ and $\,\Psi(g-h)\ge 0\,$ by Section \ref{sec.psi=} above, and then by linearity of $\,\Psi$ : $\,Id = \beta\, \Psi(g) \ge \beta\, \Psi(h) = \psi_\beta\,$. $\;\diamond$ 

\section{Proof of Proposition \ref{pro.estim}} \label{sec.compas} 

   For $\,d=1\,$ (and any $\,\eta$), it is immediate from (\ref{eqndiff3}) that \  ${\ds \psi_\beta(r) = \beta \int_0^r\left[\int_\rho^\ii\! e^{- \beta\, G}\, h\right] e^{\beta\, G(\rho)}\, d\rho}\,$, so that  in the DH case, the limit expresses as :
$$\Sigma_\beta^2 = 2 \, \b  \left( \int_{\mathbb R_+}\!\frac{e^{\b \left( 1-\sqrt{1+x^2} \right)}} {\sqrt{1+x^2}}\, dx\right)^{-1}\!\times \int_{\mathbb R_+}\! \Bigg[\int_{x}^{\infty} \frac{y \, e^{\b \left( 1-\sqrt{1+y^2} \right)}}{1+y^2}\,dy \Bigg]^2\! e^{\b \left( \sqrt{1+x^2}\,-1 \right)}\, dx = J_{\b}^{-1} K_{\b}\, , $$ 
where 
$$ J_{\b} :=\int_{\mathbb R_+}\! \frac{e^{ \b \left( 1-\sqrt{1+z^2} \right)}}{\sqrt{1+z^2}}\,dz\, , \;\, K_{\b} := 2 \b \! \int_{\mathbb R_+}\! I_{\b}(x)^2 e^{\b \left( \sqrt{1+x^2}\,-1 \right)} dx\, , \; I_{\b}(x) := \int_{x}^{\infty} \frac{y \, e^{\b \left( 1-\sqrt{1+y^2} \right)}}{1+y^2}\, dy\, . $$ 

\subsection{Behaviour as $\,\b \to \infty$} 

   Integrating by parts yields :
$$ I_{\b}(x)= \frac{e^{-\b \left( \sqrt{1+x^2}\, -1 \right)}}{\b \sqrt{1+x^2}}-\frac{1}{\b} \int_{x}^{\infty} \frac{y \times e^{- \b \left( \sqrt{1+y^2}\, -1 \right)}}{(1+y^2)^{3/2}}\, dy\, . $$
As 
$$ \left(\frac{e^{\b \left( 1-\sqrt{1+x^2}\right)}}{1+x^2}\right)^{-1}\int_{x}^{\infty} \frac{y \, e^{ \b \left( 1-\sqrt{1+y^2} \right)}}{(1+y^2)^{3/2}}\, dy\, = \int_{x}^{\infty} e^{- \b \left( \sqrt{1+y^2} -\sqrt{1+x^2} \right)} \, \frac{y \, (1+x^2)}{(1+y^2)^{3/2}}\, dy $$
$$ \leq \int_{x}^{\infty}  e^{- \b \left(\sqrt{1+y^2}-\sqrt{1+x^2}\right)} \frac{y}{\sqrt{1+y^2}}\, dy\,  =\, e^{ \b \sqrt{1+x^2}} \int_{x}^{\infty}  e^{- \b \sqrt{1+y^2}} \frac{y}{\sqrt{1+y^2}}\, dy = \frac{1}{\b}\,,$$
we get 
$$ I_{\b}(x)^2 e^{\b \left( \sqrt{1+x^2}\, -1 \right)} =\, \frac{e^{\b \left( 1- \sqrt{1+x^2}\right)}}{\b^2 \,(1+x^2)}\times \left[ 1+{\mathcal O}\left(1/\b\right)\right] . $$
Hence, 
$$ K_{\b} = 2 \, \b \, \int_{\mathbb R_+} I_{\b}(x)^2 e^{\b \left( \sqrt{1+x^2}\,-1\right)} dx \,=\,
\frac{2}{\b} \int_{\mathbb R_+} \frac{e^{\b \left( 1-\sqrt{1+x^2}\right)}}{1+x^2}\, dx
\times \left[ 1+{\mathcal O}\left(1/\b\right)\right] . $$
\noindent
Setting \  $u=\b \left( \sqrt{1+x^2}\, -1\right)$, we get : 
$$\frac{1}{\b} \int_{\mathbb R_+} \frac{e^{\b \left( 1-\sqrt{1+x^2}\right)}}{1+x^2}\, dx\, =
 \int_{\mathbb R_+} \frac{e^{-u}\, du}{(u+\b) \sqrt{u(u+2\b)}}\, \sim\,  \sqrt{\frac{\pi}{2}}\,\, \b^{-3/2} $$
by dominated convergence, as $\,\b  \to \infty\,$. 
We have similarly : 
$$ J_{\b} = \int_{\mathbb R_+}\frac{e^{\b \left( 1-\sqrt{1+z^2}\right)}}{\sqrt{1+z^2}}\, dz\, = \int_{\mathbb R_+}\frac{e^{-u}\, du}{\sqrt{u(u+2\b)}}\, \sim\, \sqrt{\frac{\pi}{2}}\,\, \b^{-1/2} \, , $$
whence 
$$ \Sigma^2_\b = J_{\b}^{-1} K_{\b}\, \sim  \, 2/\b\, . $$

\subsection{Behaviour as $\,\b \to 0$} 

   We have 
$$ J_{\b}\, =
 \int_{\mathbb R_+} \frac{e^{-u}\, du}{\sqrt{u(u+2\b)}}\, =
2 \int_{\mathbb R_+}\frac{e^{-\b u^2}du}{\sqrt{u^2+2}} =2 \left[  \int_{0}^1 \frac{e^{-\b u^2}du}{\sqrt{u^2+2}} +   \int_{0}^{1/\sqrt{\b}} \frac{e^{-1/t^2}dt}{t \sqrt{1+2\b t^2}} \right] , $$ 
where we performed the change of variable $\,\b u^2 = 1/t^2\,$.  Hence, as $\,\b \to 0\,$, we have : 
$$ J_{\b} = 2 \left[  \int_{0}^1 \frac{e^{-\b u^2}du}{\sqrt{u^2+2}} + \int_{0}^{1} \frac{e^{-1/t^2}dt}{t \sqrt{1+2\b t^2}}+ \int_{1}^{1/\sqrt{\b}} \frac{dt}{t \sqrt{1+2\b t^2}} +
\int_{1}^{1/\sqrt{\b}} \frac{ (e^{-1/t^2}-1)\, dt}{t \sqrt{1+2\b t^2}}\right]  $$ 
$$ =\, 2 \int_{1}^{1/\sqrt{\b}} \frac{dt}{t \,\sqrt{1+2\b t^2}} + C_1 + o(1) , $$ 
for a positive constant $\, C_1\,$. Integrating by parts yields then :
$$ \int_{1}^{1/\sqrt{\b}} \frac{dt}{t\, \sqrt{1+2\b t^2}}\, = \left[\frac{\log t}{\sqrt{1+2\b t^2}} \right]_1^{1/\sqrt{\b}} +\, 2\b \int_{1}^{1/\sqrt{\b}} \frac{t \log t\, dt}{\left(1+2\b t^2 \right)^{3/2}}\, . $$
Setting $\,u = \sqrt{\b}\, t\,$, we get : 
$$ \int_{1}^{1/\sqrt{\b}} \frac{dt}{t \sqrt{1+2\b t^2}}\, =\, \frac{1}{\sqrt{3}}\, \log \left( \frac{1}{\sqrt{\b}} \right) + 2 \int_{\sqrt{\b}}^1 \frac{u \log(u / \sqrt{\b})\,du}{\left(1+2u^2 \right)^{3/2}} $$
$$  =\, \frac{1}{\sqrt{3}} \log \left( \frac{1}{\sqrt{\b}} \right)
+ 2 \log\left(\frac{1}{\sqrt{\b}}\right) \int_{\sqrt{\b}}^1 \frac{u \: du}{\left(1+2u^2 \right)^{3/2}} + 2 \int_{\sqrt{\b}}^1 \frac{u \log u\, du}{\left(1+2u^2 \right)^{3/2}}\, . $$
Now, as $\,\b \to 0\,$ we have 
$$ \int_{\sqrt{\b}}^1 \frac{u \: du}{\left(1+2u^2 \right)^{3/2}} =\left[ \frac{1}{2\sqrt{1+2u^2}} \right] _{\sqrt{\b}}^1 = \frac{1}{2}\left(\frac{1}{\sqrt{1+2\b}}-\frac{1}{\sqrt{3}}\right) \longrightarrow \, \frac{1}{2} \left( 1-\frac{1}{\sqrt{3}}\right) , $$
and 
$$ 2 \int_{\sqrt{\b}}^1 \frac{u \log u\, du} {\left(1+2u^2 \right)^{3/2}}\,\longrightarrow\,
C_2 := 2 \int_{0}^1 \frac{u \log u\, du}{\left(1+2u^2 \right)^{3/2}}\, \in \mathbb R_-\, . $$
Hence, 
$$ \int_{1}^{1/\sqrt{\b}} \frac{dt}{t\, \sqrt{1+2\b t^2}}\, =\, \log \left( 1/\sqrt{\b} \right) \times [ 1+ 
o(1)] + C_2\, , $$ 
and \qquad    ${\ds J_{\b} \sim \, 2 \, \log \left(1/{\b} \right) . }$ \qquad  Otherwise, 
$$ K_{\b} = 2 \b \! \int_{\mathbb R_+}\! I_{\b}(x)^2 e^{\b \left( \sqrt{1+x^2}\,-1 \right)}dx\, =\, 2 \b  \int_{\mathbb R_+} \left( \int_{\sqrt{1+x^2}}^{\infty} e^{\b \left( 1-u\right)}\,\frac{du}{u}\right)^2 e^{\b \left(\sqrt{1+x^2}\, -1\right)} dx $$
$$ = 2  \b \int_{1}^{\infty} \left( \int_{y}^{\infty} e^{\b \left( 1-u\right)}\,\frac{du}{u}\right)^2  \frac{y \: e^{\b \left(y-1 \right)} }{\sqrt{y^2-1}}\, dy\, =\, 2 \b\, e^{\b}  \int_{1}^{\infty} \left( \int_{y}^{\infty} e^{-\b u}\, \frac{du}{u}\right)^2  \frac{y \, e^{\b y }\, dy}{\sqrt{y^2-1}}\, . $$
Observing that 
$$ \left| \left[ \int_{y}^{\infty}\! e^{-\b u}\, \frac{du}{u}\right]^2 -\left[\int_{y}^{y/\b}\! e^{-\b u}\, \frac{du}{u} \right]^2 \right| = 2 \left[\int_{y}^{y/\b}\! e^{-\b u}\, \frac{du}{u}\right] \left[\int_{y/\b}^{\infty}\! e^{-\b u}\, \frac{du}{u}\right] + \left[\int_{y/ \b }^{\infty}\! e^{-\b u}\, \frac{du}{u}\right]^2 $$
$$ \leq\,  2 \log \left( \frac{1}{\b} \right) \times \frac{e^{-y}}{y}\, + \frac{e^{-2y}}{y^2} \, , $$
and then that 
$$ \left( \log{\frac{1}{\b}}\right)^{-1} \int_{1}^{\infty} \left[\int_{y}^{\infty}\! e^{-\b u}\,\frac{du}{u} \right]^2 \frac{y\, e^{\b y}\, dy}{\sqrt{y^2-1}} = \left( \log{\frac{1}{\b}} \right)^{-1}
\int_{1}^{\infty} \left[\int_{y}^{y/\b}\! e^{-\b u}\,\frac{du}{u} \right]^2  \frac{y \, e^{\b y}\, dy}{\sqrt{y^2-1}} + {\mathcal O} \left(1 \right) , $$
changing $\,\b y\,$ into $\,y\,$ we get : 
$$ \left( \log{\frac{1}{\b}} \right)^{-1}\, K_{\b} = 2\,e^{\b}\left( \log{\frac{1}{\b}} \right)^{-1}\! A_{\b}+ {\mathcal O} (\b),  \quad \textrm{with} \quad
A_{\b} := \int_{\b}^{\infty} \left(\int_{y}^{y/\b}e^{-u}\,\frac{du}{u}\right)^2  \frac{y\, e^{y}\, dy}{\sqrt{y^2-\b^2}}\, . $$ 
Finally, setting $\,x=\sqrt{y^2-\b^2}\,$, we have by dominated convergence : 
$$A_{\b}=  \int_0^{\infty} \left( \int_{\sqrt{x^2+\b^2}}^{\sqrt{1+x^2/\b^2}}\frac{e^{- u}du}{u}\right)^2  e^{\sqrt{x^2+\b^2}}\, dx\, \longrightarrow \, A := \int_0^{\infty} \left( \int_{x}^{\infty}\frac{e^{- u}du}{u}\right)^2  e^{x}\,dx \in\mathbb R_+^*\, , $$
whence 
$$\Sigma^2_\b = \, J_{\b}^{-1} K_{\b}\, \sim\,  A \Big/ \log (1/{\b}) \, . \;\;\diamond $$

\bigskip 
\parn 
\centerline{\bf REFERENCES} 
\par\vskip 8mm 

\vbox{ \noindent 
{\bf [BDR1]} Barbachoux C. , Debbasch F. , Rivet J.P.   \  {\it Hydrodynamic behavior of Brownian particles \par \hskip 15mm  in a position-dependent constant force-field.} \quad     J. Math. Phys., vol. 40, p.2891, 2001. }
\medskip 

\vbox{ \noindent 
{\bf [BDR2]} Barbachoux C. , Debbasch F. , Rivet J.P.   \quad {\it Covariant Kolmogorov Equation and Entropy \par \hskip 2mm  Current for the Relativistic Ornstein-Uhlenbeck Process.}  \hskip 2mm         Eur. Phys. J., vol. 19, p.37, 2001. }
\medskip 

\vbox{ \noindent 
{\bf [BDR3]} Barbachoux C. , Debbasch F. , Rivet J.P.   \quad {\it The spatially one-dimensional Relativistic \par \hskip 1mm   Ornstein-Uhlenbeck process in an arbitrary inertial frame.}   \hskip 1mm  Eur. Phys. J., vol. 23, p.487, 2001. }
\medskip 

\vbox{ \noindent 
{\bf [D]} Debbasch F.   \quad {\it A diffusion process in curved space-time.} \par 
\hskip 28mm        J. Math. Phys., vol. 45, n$^o$7, p.2744-2760, 2004. }
\medskip 

\vbox{ \noindent 
{\bf [DH1]} Dunkel  J. , Hänggi P.   \  {\it Theory of relativistic Brownian Motion : The $(1+1)$-dimensional \par  \hskip 46mm   case.}   \quad   Physical Review E, vol. 71, 016124, 2005. }
\medskip 

\vbox{ \noindent 
{\bf [DH2]} Dunkel  J. , Hänggi P.   \  {\it Theory of relativistic Brownian Motion : The $(1+3)$-dimensional \par  \hskip 46mm   case.} \quad        Ar$\chi$iv, 0505532, 2005. }
\medskip 

\vbox{ \noindent 
{\bf [DMR]} Debbasch F. , Mallick K. , Rivet J.P.   \quad {\it Relativistic Ornstein-Uhlenbeck Process.} \par  \hskip 77mm        J. Stat. Phys., vol. 88, p.945, 1997. }
\medskip 

\vbox{ \noindent 
{\bf [DR]} Debbasch F. , Rivet J.P.   \  {\it A diffusion equation from the Relativistic Ornstein-Uhlenbeck \par  \hskip 49mm      Process.} \qquad   J. Stat. Phys., vol. 90, p.1179, 1998. }
\medskip 

\vbox{ \noindent 
{\bf [F]} Franchi J. \quad  {\it Relativistic Diffusion in G\"odel's Universe.} \par 
\hskip 23mm    http://arxiv.org/abs/math.PR/0612020, 2006.  } 
\medskip 

\vbox{  \noindent 
{\bf [FLJ]} Franchi J. ,  Le Jan Y.  \  {\it Relativistic Diffusions and Schwarzschild Geometry.} \par
\hskip 46mm Comm. Pure Appl. Math., vol. LX, n$^o\,$2, p.187-251, 2007.} 
\medskip 

\vbox{ \noindent 
{\bf [K]} Kallenberg O.  \quad {\it Foundation of Modern Probability.} \quad Springer, Berlin, 2001. }
\bigskip 

\vbox{ \noindent 
{\bf [RY]} Revuz D. , Yor M.  \hskip 2mm   {\it Continuous martingales and Brownian motion.} \hskip 2mm   Springer, Berlin, 1999. }
\bigskip 

\vbox{ \noindent 
{\bf [V]} Varadhan S.R.S.   \quad {\it Diffusion problems and partial differential equations.} \par 
\hskip 35mm        Springer, Berlin, 1980. }
\medskip 

\medskip 
\centerline{--------------------------------------------------------------------------------------------}
\bigskip 

\vbox{  \noindent 
J\"urgen ANGST \   and \  Jacques FRANCHI : \quad Universit\'e Louis Pasteur, I.R.M.A., \parn 
\hskip 60mm 7 rue Ren\'e Descartes, 67084 Strasbourg cedex. FRANCE. \quad \parn 
angst@math.u-strasbg.fr ,\   franchi@math.u-strasbg.fr 
}
\medskip 
\centerline{--------------------------------------------------------------------------------------------}

\end{document}